\documentclass[11pt]{article} 
\usepackage{amssymb}
\usepackage{amsthm}

\newdimen\margegd \margegd = 1 in     
\newdimen\margehb \margehb = 1.5 in   

\advance\hsize-\margegd \advance\hsize-\margegd
\advance\vsize-\margehb \advance\vsize-\margehb

\hoffset=\margegd 
\advance\hoffset-1 in 
\voffset=\margehb
\advance\voffset-1 in

\textwidth = \hsize 
\textheight = \vsize 
\advance\textheight-24bp

\swapnumbers 
\newtheorem{theorem}{Theorem}[section]
\newtheorem{prop}[theorem]{Proposition}
\newtheorem{lem}[theorem]{Lemma}

\newtheorem{definition}[theorem]{Definition}
\theoremstyle{remark} 
\newtheorem{note}[theorem]{Note}

\newcommand\Un{{\mathbf 1}} 
\newcommand\CC{{\mathbb C}}
\newcommand\NN{{\mathbb N}} 
\newcommand\aA{{\mathbb A}}
\newcommand\QQ{{\mathbb Q}} 
\newcommand\Real{\mathrm{Re}}
\newcommand\Imag{\mathrm{Im}} 
\newcommand\cF{{\cal F}}
\newcommand\cG{{\cal G}} 
\newcommand\RR{{\mathbb R}}
\newcommand\ZZ{{\mathbb Z}} 
\newcommand\HH{{\mathbb H}}
\let\wh=\widehat
\let\la=\lambda
\let\La=\Lambda

\def\From{From} 

\begin{document}

\title{On Fourier and Zeta(s)}

\author{Jean-Fran\c cois Burnol}

\date{This is the final version for FORUM MATHEMATICUM;
October 2002\footnote{\baselineskip=10pt Jan 2003: identical
with Oct 2002 except for the spelling of reminder on page 2
and the amendment of a faulty sentence on page 38 (lines
3-4) which gave an incorrect summary of general properties
of Krein type spaces. Mar 2003: see the added footnote on
page 3.}}

\maketitle

\begin{abstract}
We study some of the interactions between the Fourier Transform
and the Riemann zeta function (and Dirichlet-Dedekind-Hecke-Tate
$L$-functions).
\tableofcontents
\end{abstract}

{\baselineskip=10pt Universit\'e Lille 1, UFR de Math\'ematiques, Cit\'e
Scientifique M2, F-59655 Villeneuve d'Ascq Cedex,
France. \texttt{burnol@agat.univ-lille1.fr} \par}

\clearpage

\section{Introduction}

\subsection{A framework for this paper}

The zeta function $\zeta(s)$ assumes in Riemann's paper quite a
number of distinct identities: it appears there as a Dirichlet
series, as an Euler product, as an integral transform, as an
Hadamard product (rather, Riemann explains how $\log\zeta(s)$ may
be written as an infinite sum involving the zeros)\dots We retain
three such identities and use them as symbolic vertices for a
triangle:
$$\begin{array}{ccccc} & &\sum_n \frac{1}{n^s}& & \\ &\diagup&
&\diagdown& \\ \prod_\rho (1 - \frac{s}{\rho})& &
\frac{\qquad}{\qquad}& &\prod_p \frac{1}{1 - p^{-s}}\\
\end{array}$$

These formulae stand for various aspects of the zeta function
which, for the purposes of this manuscript, we may tentatively
name as follows:
$$\begin{array}{ccccc} & &\mbox{summations}& & \\ &\diagup&
&\diagdown& \\ \mbox{zeros}& &\overline{\qquad}& &\mbox{primes}\\
\end{array}$$
Even a casual reading of Riemann's paper will reveal how much
Fourier analysis lies at its heart, on a par with the theory of
functions of the complex variable. Let us enhance appropriately
the triangle:
$$\begin{array}{ccccc} & &\sum_n \frac{1}{n^s}& & \\
&\hfill\diagup& &\diagdown\hfill& \\ &\mathrm{Fourier}&
&\mathrm{Fourier}& \\ &\diagup\hfill& &\hfill\diagdown& \\
\prod_\rho (1 -
\frac{s}{\rho})&\overline{\qquad}&\mathrm{Fourier}&\overline{\qquad}
&\prod_p \frac{1}{1 - p^{-s}}\\
\end{array}$$
Indeed, each of the three edges is an arena of interaction
between the Fourier Transform, in various incarnations, and the
Zeta function (and Dirichlet $L$-series, or even more general
number theoretical zeta functions.) We thus specialize to a
triangle which will be the framework of this manuscript:

$$\begin{array}{ccccc} & &\begin{array}{c}\sum_n \frac{1}{n^s}\\
\mathrm{summations}\\\end{array}& & \\ &\hfil\diagup&
&\diagdown\hfill& \\ &\begin{array}{c}\mathrm{Hilbert\ spaces}\\
\mathrm{and\ vectors}\end{array}&{\mbox{\huge
?}}&\begin{array}{c}\mathrm{Adeles}\\ \mathrm{and\
Ideles}\end{array}& \\ &\diagup\hfill& &\hfill\diagdown& \\
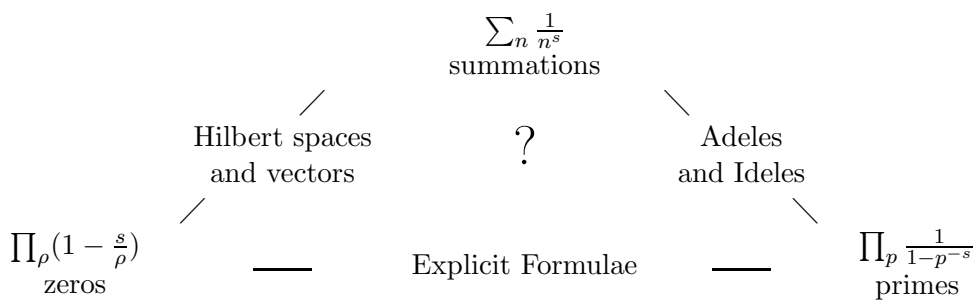
\begin{array}{c}\prod_\rho (1 - \frac{s}{\rho})\\
\mathrm{zeros}\end{array}&\overline{\qquad} &\mathrm{Explicit\
Formulae}&\overline{\qquad}&
\begin{array}{c}\prod_p \frac{1}{1 - p^{-s}}\\
\mathrm{primes}\end{array}\\
\end{array}$$

The big question mark serves as a reminder that we are missing
the $2$-cell (or $2$-cells) which would presumably be there if
the nature of the Riemann zeta function was really understood.

\subsection{The contents of this paper}

The paper contains in particular motivation, proofs, and
developments related to a ``fairly simple'' (hence
especially interesting) formula\footnote{Note added in
proofs (March 2003): the formula has in fact been discovered
earlier by Duffin and Weinberger
(\textit{Proc. Natl. Acad. Sci.}  \textbf{88} (1991),
no. 16, 7348--7350; \textit{J. Fourier Anal. Appl.} 
\textbf{3} (1997), no. 5, 487--497) and should have been
referred to here as  the ``Duffin-Weinberger dualized
Poisson formula''. Our whole analysis, which relates it to
the study of the Riemann zeta function and generalizations,
is a novel contribution.}:

$$ \int_{\RR} \Big(\sum_{n\neq 0} \frac{g(t/n)}{|n|} - \int_{\RR}
\frac{g(1/x)}{|x|}\,dx\Big)e^{2\pi i\,ut}\,dt = \sum_{m\neq 0}
\frac{g(m/u)}{|u|} - \int_{\RR} g(y)\,dy
$$

We call this the \emph{co-Poisson intertwining formula}. The
summations are over the non-zero relative integers. The formula
applies, for example, to a function $g(t)$  of class ${\cal
C}^\infty$ which is compactly supported on a closed set not
containing the origin. Then the right-hand side is a function in
the Schwartz class of smooth, rapidly decreasing functions, and
the formula exhibits it as the Fourier Transform of another
Schwartz function. These Schwartz functions have the peculiar
property of being \emph{constant, together with their Fourier
transform,} in a neighborhood of the origin. A most interesting
situation arises when the formula constructs square-integrable
functions of this type and from our discussion of this it will be
apparent that, although fairly simple, the \emph{co-Poisson
Formula} is related to a framework which is very far from being
formal.

Once found, the formula is immediately proven, and in many
different ways. Furthermore it is one among infinitely many such
co-Poisson formulae (it is planned to discuss this further in
\cite{copoisson}). This prototype is directly equivalent to the
functional equation of the Riemann zeta function. It has
implications concerning the problems of zeros.

We start with a discussion of our previous work \cite{conductor}
\cite{cras1} \cite{mrl} on the ``Explicit Formulae'' and  the
\emph{conductor operators} $\log|x|_\nu + \log|y|_\nu$. We also
include a description of our work on {adeles, ideles, scattering
and causality} \cite{imrn} \cite{jnt}, which is a first attempt
to follow from local to global the idea of multiplicatively
analysing the additive Fourier transform. This is necessary to
explain the motivation which has led to a reexamination of the
Poisson-Tate summation formula on adeles and to the discovery of
the related but subtly distinct \emph{co-Poisson intertwining}.

Both the \emph{conductor operators} and the \emph{co-Poisson
intertwining} originated from an effort to move Tate's Thesis
\cite{Tate} towards the zeros and the \emph{so-called
Hilbert-P\'olya idea}. It is notable that the zeros do not show
up at all in Tate's Thesis: the conductor operator results in
part from a continuation of the local aspects of Tate's Thesis;
the co-Poisson formula results from a reexamination of the global
aspects of Tate's Thesis. This reexamination, initially undergone
during the fall of 1998, shortly after the discovery of the
conductor operator, was also in part motivated by the preprint
version of the work of Connes \cite{Con99} (extending his earlier
Note \cite{Con96}) which had just appeared and where a very
strong emphasis is put on the \emph{so-called Hilbert-P\'olya
idea}. As we felt that the symmetries of the local conductor
operators should have some bearing on global constructions we
were very much interested by the constructions of Connes, and
especially by the attempt to realize a cut-off simultaneously in
position and momentum. This provides an indirect connection with
our work, as reported upon here. But our cut-off is (or, perhaps
better, appears to be) infrared, not ultraviolet. On our first
encounter (on the adeles) with the formula we call here
co-Poisson, we realized that we were constructing distributions
for which it was easy to compute the Fourier Transform, and that
these distributions were formally perpendicular to the zeros, but
it was not immediately apparent to us that something beyond the
usual use of the Poisson Formula was at work, as we did not at
first understand that there was a temperature parameter, and that
the Riemann zeta function is associated with a phase transition
as we vary the temperature below a certain point. So, we left
this aside for a while.

A key additional component to our effort came from the {Theorem
of B\'aez-Duarte, Balazard, Landreau and Saias} \cite{Bal3} which
is related to the Nyman-Beurling criterion \cite{Nym, Beu55} for
the validity of the Riemann Hypothesis. The link we have
established (\cite{aim}) between the so-called Hilbert-P\'olya
idea and this important {Theorem of B\'aez-Duarte, Balazard,
Landreau and Saias} leads under a further examination, which is
reported upon here, to the consideration of certain functions
which are meromorphic in the entire complex plane.

This then connects to the mechanism provided by the
\emph{co-Poisson intertwining} for the construction of {Hilbert
Spaces $H\!P_\lambda$ and Hilbert vectors $Z^\lambda_{\rho,k}$}
associated with the non-trivial zeros of the Riemann zeta
function. The method applies to Dirichlet $L$-series as well, and
the last theorem of this paper is devoted to this. Some
importance is ascribed by the author to this concluding result,
not in itself of course (as many infinitely more subtle results
than this one have been established on the zeta and $L$-functions
since Riemann's paper), but rather as a clue which could provide
inspiration for further endeavours. The light is extremely dim,
but it has the merit of  existence.

The  discussion leading to this final result makes use in
particular of an important theorem of Krein (on entire functions
of finite exponential type \cite{Kre}), and we relate the matter
with the theory of Krein type spaces as exposed in the book
\cite{Dym} by Dym and McKean. An intriguing question arises on
the properties of the \emph{Krein string} which is thus
associated with the Riemann zeta function. It seems that this
Krein string is considered here for the first time, but we add
immediately that we do not provide anything beyond mentioning it!
Rather our technical efforts, which are not
completely obvious, and not even fairly simple, lead to a
realization of the Krein type spaces of this very special Krein
string as subspaces of certain spaces $K_\lambda$
($0<\lambda<\infty$) which are involved in a kind of
multiplicative spectral (scattering) analysis of the Fourier
cosine transform. The quotient spaces $H\!P_\lambda$
($0<\lambda<1$) are the spaces we propose as an approximation
(getting better as $\lambda\to0$) to an hypothetical so-called
Hilbert-P\'olya space.

The ambient spaces $K_\lambda$ have a realization as Hilbert
Spaces of entire functions in the sense of de~Branges
\cite{Bra}. The co-Poisson formula and the discussion ot the
Nyman-Beurling criterion both suggest that it is useful to go
beyond the framework of entire functions and consider more
generally certain Hilbert spaces of meromorphic functions, but no
general development has been attempted here.

The spaces $K_\lambda$ are among the \emph{Sonine spaces}
originally studied in the sixties by de~Branges \cite{Bra64},
V.~Rovnyak \cite{Rov66} and J.~Rovnyak and V.~Rovnyak
\cite{Rov67, Rov69} (the terminology ``Sonine spaces of entire
functions'' was introduced in \cite{Rov69}). They are a special
instance of the theory of Hilbert spaces of entire functions
\cite{Bra}. But it is only for the Sonine spaces associated to
the Hankel transform of integer orders that the de~Branges
structure could be explicited in these papers. The theory of the
Sonine spaces for the cosine and sine transforms is far less
advanced. Recently though, the author has made some initial
progress on this topic (\cite{cras3}).

As de~Branges has considered the use of the general Hankel-Sonine
spaces in papers \cite{Bra92, Bra94} (and also in electronically
available unpublished manuscripts) where the matter of the
Riemann Hypothesis is mentioned, it is important to clarify that
neither the co-Poisson formula, nor the spaces $H\!P_\lambda$
($0<\lambda<1$), $W_\lambda$ and $W_\lambda^\prime$, nor the
vectors $Z^\lambda_{\rho,k}$ for $k\geq1$, have arisen in any of
de~Branges's investigations known to this author (this is said
after having spent some time to investigate the demands of the
situation created by these papers).

The circumstances of the genesis of this paper have led us to
devote a special final section, which is very brief, to some
speculations on the nature of the zeta function, the GUE
hypothesis, and the Riemann hypothesis.

\subsection{Acknowledgements}

The initial version of the manuscript, containing all essentials,
was completed in December 2001, on the occasion of the author's
``habilitation'', which took place at the University of Nice. The
author thanks Michel Balazard, Enrico Bombieri, Bernard
Candelpergher, Jean-Pierre Kahane, Philippe Maisonobe, Michel
Miniconi, and Joseph Oesterl\'e for their contribution and/or
participation. The author thanks Luis B\'aez-Duarte for
permission to incorporate a joint-proof of co-Poisson, and
Bernard Candelpergher for permission to incorporate another,
related, joint-proof of co-Poisson. The author thanks Jean-Pierre
Kahane for communicating his method of construction of Sonine
functions. The author thanks Michel Balazard and \'Eric Saias for
general discussion on the zeta function and Sonine functions.
The author thanks Michael McQuillan for support related to
matters of publication.

\section{Explicit Formulae, $\log|x|+\log|y|$, adeles,
  ideles, scattering, causality}

Riemann discovered the zeros and originated the idea of counting
the primes and prime powers (suitably weighted) using
them. Indeed this was the main focus of his famous paper. Later a
particularly elegant formula was rigorously proven by
von~Mangoldt:
$$\sum_{1<n<X}{\Lambda(n)} + \frac{1}{2}\Lambda(X) =  X -
\sum_{\rho}\frac{X^\rho}{\rho} - \log(2\pi) -
\frac{1}{2}\log(1-X^{-2})$$
\begin{small}\baselineskip=10pt
Here $X > 1$ (not necessarily an integer) and  $\Lambda(Y) =
\log(p)$ if $Y > 1$ is a positive power of the prime number $p$,
and is $0$ for all other values of $Y$. The $\rho$'s are the
Riemann Zeros (in the critical strip), the sum over them is not
absolutely convergent, even after pairing $\rho$ with
$1-\rho$. It is defined as $\lim_{T\to\infty}
\sum_{|\mathrm{Im}(\rho)| < T}{X^\rho /\rho}$.\par
\end{small}

In the early fifties Weil published a paper \cite{Wei52} on the
Riemann-von~Mangoldt Explicit Formula, and then another one
\cite{Wei72} in the early seventies which considered non-abelian
Artin (and Artin-Weil) $L$-functions. While elucidating already
in his first paper new algebraic structure, he did this
maintaining a level of generality encompassing in its scope the
von~Mangoldt formula (although it requires some steps to deduce
this formula from the Weil explicit formula.)  The analytical
difficulties arising are  an expression of the usual difficulties
with Fourier inversion. The ``test-function flavor'' of the
``Riemann-Weil explicit formula'' had been anticipated by Guinand
\cite{Gui}.

So in our opinion a more radical innovation was Weil's discovery
that the local terms of the Explicit Formulae acquire a natural
expression on the $\nu$-adics, and that this enables to put the
real and complex places on a par with the finite places (clearly
Weil was motivated by analogies with function fields, we do not
discuss that here.) Quite a lot of algebraic number theory
\cite{Wei74} is necessary in Weil's second paper to establish
this for Artin-Weil $L$-functions.

We stay here at the  simpler level of Weil's first paper and show
how to put all places of the number field at the same level. It
had first appeared in Haran's work \cite{Har} that it was
possible to formulate the Weil's local terms in a more unified
manner than had originally been done by Weil. We show that an
operator theoretical approach allows, not only to formulate, but
also to deduce the local terms in a unified manner. The starting
point is Tate's Thesis \cite{Tate}. Let $K$ be a number field and
$K_\nu$ one of its completions. Let $\chi_\nu:K_\nu^\times\to
S^1$ be a (unitary) multiplicative character. For $0<\Real(s)<1$
both $\chi_\nu(x)|x|^{s-1}$ and $\chi_\nu(x)^{-1}|x|^{-s}$ are
tempered distributions on the additive group $K_\nu$ and the
\emph{Tate's functional equations} are the identities of
distributions:
$$\cF_\nu(\chi_\nu(x)|x|^{s-1}) =
\Gamma(\chi_\nu,s)\chi_\nu(x)^{-1}|x|^{-s}$$
for certain functions $\Gamma(\chi_\nu,s)$ analytic in
$0<\Real(s)<1$, and meromorphic in the complex plane. This is the
local half of Tate's Thesis, from the point of view of
distributions. See also \cite{Gel}. Implicit in this equation is
a certain normalized choice of additive Haar measure on $K_\nu$,
and $\cF_\nu$ is the corresponding additive Fourier transform.

Let us view this from a Hilbert space perspective. The
quasi-characters $\chi_\nu(x)^{-1}|x|^{-s}$ are never
square-integrable, but for $\Real(s) = 1/2$ they are the
generalized eigenvectors arising in the spectral analysis of the
unitary group of dilations (and contractions):
$\phi(x)\mapsto\phi(x/t)/\sqrt{|t|_\nu}$, $x\in K_\nu$, $t\in
K_\nu^\times$. Let $I_\nu$ be the unitary operator
$\phi(x)\mapsto \phi(1/x)/|x|_\nu$, and let $\Gamma_\nu =
\cF_\nu\cdot I_\nu$. Then:
$$\Gamma_\nu(\chi_\nu(x)^{-1}|x|^{-s}) =
\Gamma(\chi_\nu,s)\chi_\nu(x)^{-1}|x|^{-s}$$
and this says that the $\chi_\nu(x)^{-1}|x|_\nu^{-s}$, for
$\Real(s) = 1/2$, are the generalized eigenvectors arising in the
spectral analysis of the unitary scale invariant operator
$\Gamma_\nu = \cF_\nu\cdot I_\nu$.

The question \cite{conductor} which leads from Tate's Thesis
(where the zeros do not occur at all) to the topic of the
Explicit Formulae is: \emph{what happens if we take the
derivative with respect to $s$ in Tate's functional equations?} 
Proceeding formally we obtain:
$$-\Gamma_\nu(\log|x|_\nu\; \chi_\nu(x)^{-1}|x|_\nu^{-s}) =
\Gamma^{\,\prime}(\chi_\nu,s)\cdot\chi_\nu(x)^{-1}|x|_\nu^{-s} -
\Gamma(\chi_\nu,s)\log|x|_\nu \cdot\chi_\nu(x)^{-1}|x|_\nu^{-s}$$
$$\log(|x|_\nu)\cdot\chi_\nu(x)^{-1}|x|_\nu^{-s}-\Gamma_\nu\left(\log|x|_\nu
\frac{\chi_\nu(x)^{-1}|x|_\nu^{-s}}{\Gamma(\chi_\nu,s)}\right) =
\left(\frac{d}{ds}\log\Gamma(\chi_\nu,s)\right)\cdot
\chi_\nu(x)^{-1}|x|_\nu^{-s}$$
$$\Big(\log|x|_\nu - \Gamma_\nu\cdot\log|x|_\nu\cdot
\Gamma_\nu^{-1}\Big)\cdot(\chi_\nu(x)^{-1}|x|_\nu^{-s}) =
\left(\frac{d}{ds}\log\Gamma(\chi_\nu,s)\right)\cdot
\chi_\nu(x)^{-1}|x|_\nu^{-s}$$
Let $H_\nu$ be the scale invariant operator $\log|x|_\nu -
\Gamma_\nu\cdot\log|x|_\nu\cdot \Gamma_\nu^{-1} = \log|x|_\nu +
\cF_\nu\cdot\log|x|_\nu\cdot\cF_\nu^{-1}$, which we also write
symbolically as:
$$H_\nu = \log|x|_\nu + \log|y|_\nu$$ then we see that the
conclusion is:

\begin{theorem}[\cite{conductor} \cite{cras1}]
The generalized eigenvalues of the \emph{conductor operator}
$H_\nu$ are the logarithmic derivatives of the Tate Gamma
functions:
$$H_\nu(\chi_\nu(x)^{-1}|x|_\nu^{-s}) =
\left(\frac{d}{ds}\log\Gamma(\chi_\nu,s)\right)\cdot
\chi_\nu(x)^{-1}|x|_\nu^{-s}
$$
\end{theorem}

Let $g(u)$ be a smooth function with compact support in
$\RR_+^\times$. Let $\wh{g}(s) = \int g(u)u^{s-1}\,du$ be its
Mellin transform. Let $\chi$ be a unitary character on the idele
class group of the number field $K$, with local components
$\chi_\nu$. Let $Z(g,\chi)$ be the sum of the values of
$\wh{g}(s)$ at the zeros of the (completed) Hecke $L$-function
$L(\chi,s)$ of $\chi$, minus the contribution of the poles when
$\chi$ is a principal character ($t\mapsto|t|^{-i\tau}$). Using
the calculus of residues we obtain $Z(g,\chi)$ as the integral of
$\wh{g}(s)(d/ds)\log L(\chi,s)$ around the contour of the
infinite rectangle $-1\leq \Real(s)\leq 2$. It turns out that the
compatibility between Tate's Thesis (local half) and Tate's
Thesis (global half) allows to use the functional equation
without having ever to write down explicitely all its details
(such as the discriminant of the number field and the conductor
of the character), and leads to:

\begin{theorem}[\cite{conductor}]
The explicit formula is given by the logarithmic derivatives of
the Tate Gamma functions:
$$Z(g,\chi) = \sum_\nu \int_{\Real(s)=\frac{1}{2}}
\left(\frac{d}{ds}\log\Gamma(\chi_\nu,s)\right)
\wh{g}(s)\,\frac{|ds|}{2\pi}$$
\end{theorem}

At an archimedean place the values $\wh{g}(s)$ on the critical
line give the multiplicative spectral decomposition of the
function $g_{\chi,\nu}:= x\mapsto \chi_\nu(x)^{-1}g(|x|_\nu)$ on
(the additive group) $K_\nu$, and, after checking normalization
details, one finds that the local term has exact value
$H_\nu(g_{\chi,\nu})(1)$. At a non-archimedean place, one
replaces the integral on the full critical line with an integral
on an interval of periodicity of the Tate Gamma function, and
applying Poisson summation (in the vertical direction) to
$\wh{g}(s)$ to make it periodical as well it is seen to transmute
into the multiplicative spectral decomposition of the function
$g_{\chi,\nu}:= x\mapsto \chi_\nu(x)^{-1}g(|x|_\nu)$ on
$K_\nu$\thinspace! So we jump directly from the critical line to
the completions of the number field $K$, and end up with the
following version of the explicit formula:

\begin{theorem}[\cite{conductor} \cite{cras1}]
Let at each place $\nu$ of the number field $K$:
$$g_{\chi,\nu}= x\mapsto \chi_\nu(x)^{-1}g(|x|_\nu)$$ on $K_\nu$
($g_{\chi,\nu}(0) = 0$). Then
$$Z(g,\chi) = \sum_\nu H_\nu(g_{\chi,\nu})(1)$$
 where $H_\nu$ is the scale invariant operator $\log|x|_\nu +
 \log|y|_\nu$ acting on $L^2(K_\nu, dx_\nu)$.
\end{theorem}

As we evaluate at $1$, the ``$\log|x|_\nu$'' half of $H_\nu$
could be dropped, and we could sum up the situation as follows:
\emph{Weil's local term is the (additive) Fourier transform of
the logarithm\thinspace!} This is what Haran had proved
(\cite{Har}, for the Riemann zeta function), except that he
formulated this in terms of Riesz potentials $|y|_\nu^{-s}$, and
did a separate check for finite places and the infinite place
that the Weil local terms may indeed be written in this way.
The explicit formula as stated above with the help of the
operator $\log|x|_\nu + \log|y|_\nu$ incorporates in a more
visible manner the compatibility with the functional
equations. Indeed we have

\begin{theorem}[\cite{conductor} \cite{cras1}]
The conductor operator $H_\nu$ commutes with the operator $I_\nu$:
$$H_\nu\cdot I_\nu = I_\nu\cdot H_\nu$$
or equivalently as $I_\nu\cdot\log|x|_\nu\cdot I_\nu =
-\log|x|_\nu$:
$$I_\nu\cdot\log|y|_\nu\cdot I_\nu = 2\log|x|_\nu + \log|y|_\nu$$
\end{theorem}

To see abstractly why this has to be true, one way is to observe
that $H_\nu$ and $\Gamma_\nu$ are simultaneously diagonalized by
the multiplicative characters, hence they commute. But obviously
$H_\nu$ commutes with $\cF_\nu$ so it has to commute with
$I_\nu$. Later, when dealing with what we call ``co-Poisson
intertwining'', we will see a similar argument in another context.

Let us suppose $\chi_\nu$ to be ramified (which means not trivial
when restricted to the $\nu$-adic units) and let $f(\chi_\nu)$ be
its \emph{conductor exponent}, $e_\nu$ the number field
\emph{differental exponent} at $\nu$, and $q_\nu$  the
cardinality of the residue field. In Tate's Thesis \cite{Tate},
one finds for a ramified character
$$\Gamma(\chi_\nu,s) = w(\chi_\nu) q_\nu^{(f(\chi_\nu) +
e_\nu)(s-\frac12)}$$
where $w(\chi_\nu)$ is a certain non-vanishing complex number,
quite important in Algebraic Number Theory, but not here. Indeed
we take the logarithmic derivative and find:
$$\frac{d}{ds}\log\Gamma(\chi_\nu,s) = (f(\chi_\nu) + e_\nu)\log
q_\nu$$
So that:
$$H_\nu(\chi_\nu^{-1}(x)\Un_{|x|_\nu=1}(x)) =  (f(\chi_\nu) +
e_\nu)\log(q_\nu)\chi_\nu^{-1}(x)\Un_{|x|_\nu=1}(x)$$
 which says that ramified characters are eigenvectors of $H_\nu$
 with eigenvalues $(f(\chi_\nu) + e_\nu)\log q_\nu$. Hence the
 name ``conductor operator'' for $H_\nu$. We note that this
 contribution of the differental exponent is there also for a
 non-ramified character and explains why in our version of the
 Explicit  Formula there is no explicit presence of the
 discriminant of the number field. If we now go through the
 computation of the distribution theoretic additive Fourier
 transform of $\log|x|_\nu$ and compare with the above we end up
 with a proof \cite{conductor} of the well-known Weil integral
 formula \cite{Wei52} \cite{Wei72} \cite{Wei74} (Weil writes
 $d^\times t$ for $\log(q_\nu)d^*t$):
$$f(\chi_\nu)\log q_\nu  = \int_{K_\nu^\times} \Un_{|t|_\nu=1}(t)
\frac{1 - \chi_\nu(t)}{|1 - t|_\nu}\,d^\times t$$

In Weil's paper \cite{Wei52} we see that this formula's r\^ole
has been somewhat understated. Clearly it was very important to
Weil as it confirmed that it was possible to express similarly
all contributions to the Explicit Formula: from the infinite
places, from finite unramified places, and from finite ramified
places. Weil leaves establishing the formula to the attentive
reader. In his second paper \cite{Wei72} he goes on to extend the
scope to Artin $L$-function, and this is far from an obvious
thing.

Let us now consider the zeta and $L$-functions from the point of
view of Adeles and Ideles.   Again a major input is Tate's
Thesis. There the functional equations of the abelian
$L$-functions are established in a unified manner, but the zeros
do not appear at all. It is only recently that progress on this
arose, in the work of Connes \cite{Con99}. We have examined this
question anew \cite{imrn} \cite{jnt}, from the point of view of
the study of the interaction between the additive and
multiplicative Fourier Transforms \cite{conductor} \cite{imrn},
which as we saw is a mechanism underlying the
operator theoretic approach to the explicit formula. This led us
to the scattering theory of Lax and Phillips \cite{Lax67} and to
a formulation of the Riemann Hypothesis, simultaneous for all
$L$-functions, as a property of \emph{causality} \cite{jnt}

A key theorem from the global half of Tate's Thesis is the
following:
$$\sum_{q\in K} \cF(\varphi)(q v) = \frac{1}{|v|}\sum_{q\in K}
\varphi(\frac{q}{v})$$
This was called the ``Riemann-Roch Theorem'' by Tate, but we prefer to
call it the \emph{Poisson-Tate formula} (which sounds less definitive,
and more to the point). Let us explain the notations: $K$ is a number
field, $\varphi(x)$ is a function on the {adeles} $\aA$ of $K$
(satisfying suitable conditions), $q\in K$ is diagonally considered as
an adele, $v\in \aA^\times$ is an {idele} and $|v|$ is its
module. Finally $\cF$ is the adelic additive Fourier Transform (we
refer to \cite{Tate} for the details of the
normalizations\footnote{e.g., on $\RR$ the Fourier transform is
$\widetilde\varphi(y) = \int_\RR e^{2\pi
iyx}\varphi(x)\,dx$}\setcounter{footnote}0). We note that it does not
matter if we exchange the $\varphi$ on the right with the
$\cF(\varphi)$ on the left as $\cF(\cF(\varphi))(x) = \varphi(-x)$ and
$-1\in K^\times$. A suitable class of functions stable under $\cF$ for
which this works is given by the Bruhat-Schwartz functions: finite
linear combinations of infinite product of local factors, almost all
of them being the indicator function of the local ring of integers, in
the Schwartz class for the infinite places, locally constant with
compact support at each finite place. To each such function and
unitary character $\chi$ on the idele class group Tate associates an
$L$-function $L(\chi,\varphi)(s) = \int_{\mathrm{ideles}}
\varphi(v)\chi(v)|v|^s d^*v$, and shows how to choose $\varphi$ so
that this coincides exactly with the complete Hecke $L$-function with
grossencharakter $\chi$.

Let us write $E_0$ (very soon we will switch to a related $E$)
for the map which to the function $\varphi(x)$ on the
\emph{adeles} associates the function $\sum_{q\in K} \varphi(q
v)\sqrt{|v|}$ on the \emph{ideles} or even on the \emph{idele
class group} ${\cal C}_K$ (ideles quotiented by $K^\times$). This
map $E_0$ plays an important r\^ole in the papers of Connes
\cite{Con96} \cite{Con99} (where it is used under the additional
assumption $\varphi(0) = 0 = \cF(\varphi)(0)$, and then coincides
with the $E$ we introduce next.)  The Poisson-Tate formula tells
us that $E_0$ intertwines the additive Fourier transform with the
operator $I:g(v)\mapsto g(1/v)$.
$$(E_0\cdot\cF)(\varphi)= (I\cdot E_0)(\varphi)$$

If we (\cite{jnt}) manipulate a little bit the Poisson-Tate
formula into:
$$\sqrt{|v|}\sum_{q\in K^\times} \cF(\varphi)(q v) -
\frac{1}{\sqrt{|v|}}\varphi(0) = \frac{1}{\sqrt{|v|}}\sum_{q\in
K^\times} \varphi(\frac{q}{v}) - \sqrt{|v|}\int_\mathrm{adeles}
\varphi(x)\,dx$$
 we still have the intertwining property
$$(E\cdot\cF)(\varphi)= (I\cdot E)(\varphi)$$
 where we have  written $E$ for the map which to $\varphi(x)$
 associates 
$$u\mapsto \sqrt{|v|}\sum_{q\in K^\times} \varphi(q v)
 - \frac{\int_\mathrm{adeles}\varphi(x)dx}{\sqrt{|u|}}$$
 on the idele class group ($v$ in the class $u$).

\begin{note}
Let $\phi(x)$ be an even Schwartz function on $\RR$. Let $F(x) =
\sum_{n\geq1} \phi(nx)$. We have $\lim_{x\to0} |x|F(x) =
\int_0^\infty \phi(x)\,dx$. For $\Real(s)>1$ we may intervert the
integral with the summation and this gives $\int_0^\infty
F(x)x^{s-1}\,dx = \zeta(s)\int_0^\infty \phi(x)x^{s-1}\,dx$. The
analytic continuation of this formula to the critical strip
($0<\Real(s)<1$) requires a modification which is due to M\"untz
(as stated in Titchmarsh's book \cite[II.11]{Titchmarsh})
$$\int_0^\infty  \Big(F(x) - \frac{\int_0^\infty
\phi(y)\,dy}{x}\Big)\,x^{s-1}\,dx = \zeta(s)\int_0^\infty
\phi(x)x^{s-1}\,dx$$
So it is in truth not the original Poisson summation but the
M\"untz-modified Poisson (where one takes out $\phi(0)$ and
replaces it with $-(\int_\RR \phi(y)\,dy)/|x|)$) which
corresponds to $\zeta(s)$ as multiplier. The M\"untz modification
was used by the author in \cite{jnt} unknowingly of its previous
appearance in the literature. The author thanks Luis
B\'aez-Duarte for pointing out the reference to the section of
the book of Titchmarsh where the M\"untz formula is discussed.
\end{note}

\begin{theorem}[\cite{jnt}]
For $\varphi$ a Bruhat-Schwartz function $E(\varphi)$ is
square-integrable on the idele class group ${\cal C}_K$ (for the
multiplicative Haar measure $d^*u$), and its unitary
multiplicative Fourier transform, as a function of the unitary
characters, coincides up to an overall constant with the Tate
$L$-function on the critical line $\Real(s) = 1/2$. The functions
$E(\varphi)$ are dense in $L^2({\cal C}_K, d^*u)$ and
$E(\cF(\varphi)) = I(E(\varphi))$.
\end{theorem}

Connes \cite{Con96} \cite{Con99} had already considered the
functions $\sqrt{|v|}\sum_{q\in K^\times} \varphi(q v)$ with
$\varphi(0) = 0 = \cF(\varphi)(0)$ and he had shown that they are
dense in $L^2({\cal C}_K, d^*u)$. Let ${\cal S}_1$ be the set of
Bruhat-Schwartz functions $\varphi(x)$ which are supported in a
parallelepiped $P(v) = \{\forall\nu\;|x|_\nu\leq|v|_\nu\}$ with
$|v|\leq 1$. Let ${\cal D}_+ = E({\cal S}_1)^\perp$ and let
${\cal D}_- = E(\cF({\cal S}_1))^\perp$. The following holds:

\begin{theorem}[\cite{jnt}]
The subspaces ${\cal D}_+$ and ${\cal D}_-$ of the Hilbert space
of square-integrable functions on the idele class group are
outgoing and incoming subspaces for a Lax-Phillips scattering
system, where the idele class group plays the r\^ole of time. The
Riemann Hypothesis for all abelian $L$-functions of $K$ holds if
and only if the causality axiom ${\cal D}_+\perp{\cal D}_-$ is
satisfied.
\end{theorem}

\section{Poisson-Tate and a novel relative: co-Poisson}

We have already mentioned the Tate $L$-functions (the integral is
absolutely convergent for $\Real(s)>1$):
$$L(\chi,\varphi)(s) = \int_{\mathrm{ideles}}
\varphi(v)\chi(v)|v|^s d^*v$$

Using the Poisson-Tate summation formula, Tate established the
analytic continuation and the functional equations:
$$L(\chi, \cF(\varphi))(s) = L(\chi^{-1}, \varphi)(1-s)$$
This follows from an integral representation
$$L(\chi,\varphi)(s) = C\delta_\chi\left(
\frac{\cF(\varphi)(0)}{s-1-i\tau} -
\frac{\varphi(0)}{s-i\tau}\right) + \int_{|v|\geq1}
\Big(\varphi(v)\chi(v)|v|^s +
\cF(\varphi)(v)\chi(v)^{-1}|v|^{1-s}\Big) \;d^*v$$
 where $C$ is a certain constant associated to the number field
 $K$ (and relating the Haar measures $d^*v$ on $\aA^\times$ and
 $d^*u$ on ${\cal C}_K$), and the Kronecker symbol $\delta_\chi$
 is $1$ or $0$ according to whether $\chi(v) = |v|^{-i\tau}$ for
 a certain $\tau\in\RR$ (principal unitary character) or not
 (ramified unitary character). The integral over the ideles (this
 is not an integral over the idele \emph{classes}) with
 $|v|\geq1$ is absolutely convergent for \emph{all} $s\in\CC$.

We recall that we  associated to the Bruhat-Schwartz function
$\varphi(x)$ on the adeles  the square-integrable function
$E(\varphi)(u)$ on the idele class group ${\cal C}_K$ (with
$u\in{\cal C}_K$ the class of $v\in \aA^\times$):
$$E(\varphi)(u) = \sqrt{|u|}\sum_{q\in K^\times} \varphi(q v) -
\frac{\int_\mathrm{adeles}\varphi(x)\,dx}{\sqrt{|u|}}$$
The precise relation \cite{jnt} to the Tate $L$-functions is:
$$L(\chi,\varphi)(s) = C\int_{{\cal C}_K}
E(\varphi)(u)\chi(u)|u|^{s-1/2}\,d^*u$$
 This integral representation is absolutely convergent for
 $0<\Real(s)<1$ and we read the functional equations directly
 from it and from the intertwining property $E\cdot\cF = I\cdot
 E$.

Let $\Delta(u)$ be the  function of $u$ with values in the
distributions on the adeles:
$$\Delta(u)(\varphi) = C\cdot\left(\sqrt{|u|}\sum_{q\in K^\times}
\varphi(q v) -
\frac{\int_\mathrm{adeles}\varphi(x)\,dx}{\sqrt{|u|}}\right) =
C\cdot E(\varphi)(u)$$

We will show that it is relevant to look at $\Delta(u)$ not as a
\emph{function} in $u$ (which it is from the formula above) but
as a \emph{distribution} in $u$ (so that $\Delta$ is a
distribution with values in distributions\dots) It takes  time to
explain why this is not a tautological change of
perspective. Basically we shift the emphasis from the
Poisson-Tate summation \cite{Tate} \cite{Con99} \cite{jnt}  which
goes from adeles to ideles, to the \emph{co-Poisson summation}
which goes from ideles to adeles. The Poisson-Tate summation is a
function with values in distributions, whereas the
co-Poisson-Tate summation is a distribution whose values we try
to represent as $L^2$-functions.

Let $g(v)$ be a compact Bruhat-Schwartz function on the idele
group $\aA^\times$. This is a finite linear combination of
infinite products $v\mapsto\prod_\nu g_\nu(v_\nu)$, where almost
each component is the indicator function of the  $\nu$-adic
units, the component $g_\nu(v_\nu)$ at an infinite place is a
smooth compactly supported function on $K_\nu^\times$, and the
components at finite places are locally constant compactly
supported.

\begin{definition}
The \emph{co-Poisson summation} is the map $E^\prime$ which
assigns to each compact Bruhat-Schwartz function $g(v)$ the
distribution on the \emph{adeles} given by:
$$E^\prime(g)(\varphi) = \int_{\aA^\times} \varphi(v)\sum_{q\in
K^\times} g(qv) \sqrt{|v|}\,d^*v - \int_{\aA^\times}
g(v)|v|^{-1/2}d^*v\; \int_{\aA} \varphi(x)\,dx$$
\end{definition}

\begin{note}
Clearly $E^\prime(g)$ depends on $g(v)$ only through the function
$R(g)$ on ${\cal C}_K$ given by $R(g)(u) = \sum_{q\in K^\times}
g(qv)$, with $u\in{\cal C}_K$ the class of $v\in
\aA^\times$. However, for various reasons (among them avoiding
the annoying constant  $C$ in all our formulae), it is better to
keep the flexibility provided by $g$. The function $R(g)$ has
compact support. To illustrate this with an example, and explain
why the integral above makes sense, we take $K=\QQ$, $g(v) =
\prod_p \Un_{|v_p|_p=1}(v_p)\cdot g_\infty(v_\infty)$. Then
$\sum_{q\in \QQ^\times} g(qv) = g_\infty(|v|) + g_\infty(-|v|)$.
We may bound this from above by  a multiple of $|v|$ (as
$g_\infty$ has compact support in $\RR^\times$), and the integral
$\int_{\aA^\times} \varphi(v)|v|^{3/2}\,d^*v$ converges
absolutely as $1<3/2$ (we may take $\varphi(x)$ itself to be an
infinite product here.)
\end{note}

\begin{theorem}
The co-Poisson summation intertwines the operator $I:g(v)\mapsto
g(1/v)$ with the additive adelic Fourier Transform $\cF$:
$$\cF(E^\prime(g)) = E^\prime(I(g))$$
 Furthermore it intertwines between the multiplicative
 translations $g(v)\mapsto g(v/w)$ on ideles and the
 multiplicative translations on adelic distributions $D(x)\mapsto
 D(x/w)/\sqrt{|w|}$. And the distribution $E^\prime(g)$ is
 invariant under the action of the multiplicative group
 $K^\times$ on the adeles.
\end{theorem}

\begin{proof}
We have:
\begin{eqnarray*}
E^\prime(g)(\varphi)&=& C\int_{{\cal C}_K} \sum_{q\in K^\times}
\varphi(qv) R(g)(u)\sqrt{|u|}\,d^*u - \int_{\aA^\times}
g(v)|v|^{-1/2}d^*v\int_{\aA} \varphi(x)\,dx\\ &=& C\int_{{\cal
C}_K} \left(\frac{E(\varphi)(u)}{\sqrt{|u|}} + \frac{\int
\varphi(x)\,dx}{|u|}\right) R(g)(u)\sqrt{|u|}\,d^*u -
\int_{\aA^\times}\frac{g(v)}{\sqrt{|v|}}\,d^*v \int_{\aA}
\varphi(x)\,dx\\ &=& C\int_{{\cal C}_K}
E(\varphi)(u)R(g)(u)\,d^*u + \Big(C\int_{{\cal
C}_K}\frac{R(g)(u)}{\sqrt{|u|}}\,d^*u - \int_{\aA^\times}
\frac{g(v)}{\sqrt{|v|}}d^*v \Big)\int \varphi(x)\,dx\\ &=&
C\int_{{\cal C}_K} E(\varphi)(u)R(g)(u)\,d^*u\\
\end{eqnarray*}
Using the intertwinings $E\cdot \cF = I\cdot E$ and  $R(I(g))(u)
= R(g)(1/u)$ we get:
$$E^\prime(g)(\cF(\varphi)) = C\int_{{\cal C}_K}
E(\varphi)(u)R(g)(\frac{1}{u})\,d^*u =E^\prime(I(g))(\varphi)$$
 which completes the proof of $\cF(E^\prime(g)) =
 E^\prime(I(g))$. The compatibility with multiplicative
 translations is easy, and the invariance under the
 multiplication $x\mapsto qx$ follows.
\end{proof}

\begin{note}
The way the ideles have to act on the distributions on adeles for
the intertwining suggests some Hilbert space properties of the
distribution $E^\prime(g)$ (more on this later).
\end{note}

\begin{note}
The invariance under $K^\times$ suggests that it could perhaps be
profitable to discuss $E^\prime(g)$ from the point of view of the
Connes space $\aA/K^\times$ \cite{Con96} \cite{Con99}.
\end{note}

\begin{theorem}
The following Riemann-Tate formula holds:
$$E^\prime(g)(\varphi) = \int_{|v|\geq1}
\cF(\varphi)(v)\sum_{q\in K^\times} g(q/v)\sqrt{|v|}d^*v +
\int_{|v|\geq1} \varphi(v)\sum_{q\in K^\times}
g(qv)\sqrt{|v|}d^*v$$
$$-\varphi(0)\int_{|v|\geq1} g(1/v)|v|^{-1/2}d^*v -\int_\aA
\varphi(x)\,dx\cdot\int_{|v|\geq1} g(v)|v|^{-1/2}d^*v$$
\end{theorem}

\begin{proof}
\From\  $E^\prime(g)(\varphi) = C\int_{{\cal C}_K}
E(\varphi)(u)R(g)(u)\,d^*u$ we get
\begin{eqnarray*}
E^\prime(g)(\varphi) &=& C\int_{|u|\leq1}
E(\varphi)(u)R(g)(u)\,d^*u + C\int_{|u|\geq1}
E(\varphi)(u)R(g)(u)\,d^*u\\ &=&C\int_{|u|\geq1}
E(\cF(\varphi))(u)R(g)(1/u)\,d^*u + C\int_{|u|\geq1}
E(\varphi)(u)R(g)(u)\,d^*u\\ &=
&C\int_{|u|\geq1}\Big(E(\cF(\varphi))(u) +
\frac{\varphi(0)}{\sqrt{|u|}}\Big)R(g)(1/u)\,d^*u
-\varphi(0)\int_{|v|\geq1} g(1/v)|v|^{-1/2}d^*v\\ &+&
C\int_{|u|\geq1}\Big(E(\varphi)(u) + \frac{\int_\aA
\varphi(x)\,dx}{\sqrt{|u|}}\Big)R(g)(u)\,d^*u -(\int
\varphi)\int_{|v|\geq1} g(v)|v|^{-1/2}d^*v \\ &=& \int_{|v|\geq1}
\cF(\varphi)(v)\sum_{q\in K^\times} g(q/v)\sqrt{|v|}d^*v +
\int_{|v|\geq1} \varphi(v)\sum_{q\in K^\times}
g(qv)\sqrt{|v|}d^*v\\ &-&\varphi(0)\int_{|v|\geq1}
g(1/v)|v|^{-1/2}d^*v  -(\int_\aA \varphi(x)\,dx)\int_{|v|\geq1}
g(v)|v|^{-1/2}d^*v \\
\end{eqnarray*}
which completes the proof.
\end{proof}

We note that if we replace formally $\sum_{q\in K^\times}
g(qv)\sqrt{|v|}$ with $\chi(v)|v|^s$ we obtain exactly the Tate
formula for $L(\chi,\varphi)(s)$. But some new flexibility arises
with a ``compact'' $g(v)$:

\begin{theorem}
The following formula holds:
$$E^\prime(g)(\varphi) = \int_{|v|\leq1}
\cF(\varphi)(v)\sum_{q\in K^\times} g(q/v)\sqrt{|v|}d^*v +
\int_{|v|\leq1} \varphi(v)\sum_{q\in K^\times}
g(qv)\sqrt{|v|}d^*v -$$
$$-\varphi(0)\int_{|v|\leq1} g(1/v)|v|^{-1/2}d^*v -\int_\aA
\varphi(x)\,dx\cdot\int_{|v|\leq1} g(v)|v|^{-1/2}d^*v$$
\end{theorem}

\begin{proof}
Exactly the same as above exchanging everywhere $|u|\geq1$ with
$|u|\leq1$ (we recall that $R(g)$ has compact support and also
that as we are dealing with a number field $|v|=1$ has zero
measure). This is not possible with a quasicharacter in the place
of $R(g)(u)$. Alternatively one adds to the previous formula and
checks that one obtains $2E^\prime(g)(\varphi)$ (using
$E^\prime(I(g))(\cF(\varphi)) = E^\prime(g)(\varphi)$).
\end{proof}

To illustrate some Hilbert Space properties of the co-Poisson
summation, we will assume $K=\QQ$. The components $(a_\nu)$ of an
adele $a$ are written $a_p$ at finite places and $a_r$ at the
real place. We have an embedding of the Schwartz space of
test-functions on $\RR$ into the Bruhat-Schwartz space on $\aA$
which sends $\psi(x)$ to $\varphi(a) = \prod_p
\Un_{|a_p|_p\leq1}(a_p)\cdot \psi(a_r)$, and we write
$E^\prime_\RR(g)$ for the distribution on $\RR$ thus obtained
from $E^\prime(g)$ on $\aA$.

\begin{theorem}
Let $g$ be a compact Bruhat-Schwartz function on the ideles of
$\QQ$. The co-Poisson summation $E^\prime_\RR(g)$ is a
square-integrable function (with respect to the Lebesgue
measure). The $L^2(\RR)$ function $E^\prime_\RR(g)$ is equal to
the constant $-\int_{\aA^\times} g(v)|v|^{-1/2}d^*v$ in a
neighborhood of the origin.
\end{theorem}

\begin{proof}
We may first, without changing anything to $E^\prime_\RR(g)$,
replace $g$ with its average under the action of the finite unit
ideles, so that it may be assumed invariant. Any such compact
invariant $g$ is a finite linear combination of suitable
multiplicative translates of functions of the type $g(v) =
\prod_p \Un_{|v_p|_p=1}(v_p)\cdot f(v_r)$ with $f(t)$ a smooth
compactly supported function on $\RR^\times$, so that we may
assume that $g$ has this form.   We claim that:
$$\int_{\aA^\times} |\varphi(v)|\sum_{q\in \QQ^\times} |g(qv)|
\sqrt{|v|}\,d^*v < \infty$$
 Indeed $\sum_{q\in \QQ^\times} |g(qv)| = |f(|v|)| +|f(-|v|)|$ is
 bounded above by a multiple of $|v|$. And $\int_{\aA^\times}
 |\varphi(v)||v|^{3/2}\,d^*v < \infty$ for each Bruhat-Schwartz
 function on the adeles (basically, from $\prod_p (1 -
 p^{-3/2})^{-1} < \infty$). So
$$E^\prime(g)(\varphi) = \sum_{q\in \QQ^\times} \int_{\aA^\times}
\varphi(v)g(qv) \sqrt{|v|}\,d^*v - \int_{\aA^\times}
\frac{g(v)}{\sqrt{|v|}}d^*v\; \int_{\aA} \varphi(x)\,dx$$
$$E^\prime(g)(\varphi) = \sum_{q\in \QQ^\times} \int_{\aA^\times}
\varphi(v/q)g(v) \sqrt{|v|}\,d^*v - \int_{\aA^\times}
\frac{g(v)}{\sqrt{|v|}}d^*v\; \int_{\aA} \varphi(x)\,dx$$
Let us now specialize to $\varphi(a) = \prod_p
\Un_{|a_p|_p\leq1}(a_p)\cdot \psi(a_r)$. Each integral can be
evaluated as an infinite product. The finite places contribute
$0$ or $1$ according to whether $q\in\QQ^\times$ satisfies
$|q|_p<1$ or not. So only the inverse integers $q = 1/n$,
$n\in\ZZ$, contribute:
$$E^\prime_\RR(g)(\psi) = \sum_{n\in\ZZ^\times} \int_{\RR^\times}
\psi(nt)f(t) \sqrt{|t|}\frac{dt}{2|t|} - \int_{\RR^\times}
\frac{f(t)}{\sqrt{|t|}}\frac{dt}{2|t|}\; \int_{\RR} \psi(x)\,dx$$
We can now revert the steps, but this time on $\RR^\times$ and we
get:
$$E^\prime_\RR(g)(\psi) = \int_{\RR^\times}
\psi(t)\sum_{n\in\ZZ^\times} \frac{f(t/n)}{\sqrt{|n|}}
\frac{dt}{2\sqrt{|t|}} - \int_{\RR^\times}
\frac{f(t)}{\sqrt{|t|}}\frac{dt}{2|t|}\; \int_{\RR} \psi(x)\,dx$$
 Let us express this in terms of $\alpha(y) = (f(y)+
 f(-y))/2\sqrt{|y|}$:
$$E^\prime_\RR(g)(\psi) = \int_\RR \psi(y)\sum_{n\geq1}
\frac{\alpha(y/n)}{n}dy - \int_0^\infty \frac{\alpha(y)}{y}dy\;
\int_{\RR} \psi(x)\,dx$$
 So the distribution $E^\prime_\RR(g)$ is in fact the even smooth
 function
$$E^\prime_\RR(g)(y) = \sum_{n\geq1} \frac{\alpha(y/n)}{n} -
\int_0^\infty \frac{\alpha(y)}{y}dy$$
 As $\alpha(y)$ has compact support in $\RR\setminus\{0\}$, the
 summation over $n\geq1$ contains only vanishing terms for $|y|$
 small enough. So $E^\prime_\RR(g)$ is equal to the constant $-
 \int_0^\infty \frac{\alpha(y)}{y}dy = - \int_{\RR^\times}
 \frac{f(y)}{\sqrt{|y|}}\frac{dy}{2|y|} = -\int_{\aA^\times}
 g(t)/\sqrt{|t|}\,d^*t$ in a neighborhood of $0$. To prove that
 it is $L^2$, let $\beta(y)$ be the smooth compactly supported
 function $\alpha(1/y)/2|y|$ of $y\in\RR$ ($\beta(0) = 0$). Then
 ($y\neq0$):
$$E^\prime_\RR(g)(y) = \sum_{n\in\ZZ}
\frac{1}{|y|}\beta(\frac{n}{y}) - \int_\RR \beta(y)\,dy$$
\From\  the usual Poisson summation formula, this is also:
$$\sum_{n\in\ZZ} \gamma(ny) - \int_\RR \beta(y)\,dy = \sum_{n\neq
0} \gamma(ny)$$
 where $\gamma(y) = \int_\RR \exp(i\,2\pi yw)\beta(w)\,dw$ is a
 Schwartz rapidly decreasing function. From this formula we
 deduce easily that $E^\prime_\RR(g)(y)$ is itself in the
 Schwartz class of rapidly decreasing functions, and in
 particular it is is square-integrable.
\end{proof}

It is useful to recapitulate some of the results arising in this
proof:

\begin{theorem}\label{copoisson}
Let $g$ be a compact Bruhat-Schwartz function on the ideles of
$\QQ$. The co-Poisson summation $E^\prime_\RR(g)$ is an even
function on $\RR$ in the Schwartz class of rapidly decreasing
functions. It is constant, as well as its Fourier Transform, in a
neighborhood of the origin. It may be written as
$$E^\prime_\RR(g)(y) = \sum_{n\geq1} \frac{\alpha(y/n)}{n} -
\int_0^\infty \frac{\alpha(y)}{y}dy$$
 with a function $\alpha(y)$ smooth with compact support away
 from the origin, and conversely each such formula corresponds to
 the co-Poisson summation $E^\prime_\RR(g)$  of a compact
 Bruhat-Schwartz function on the ideles of $\QQ$. The Fourier
 transform $\int_\RR E^\prime_\RR(g)(y)\exp(i2\pi wy)\,dy$
 corresponds in the formula above to the replacement
 $\alpha(y)\mapsto \alpha(1/y)/|y|$.
\end{theorem}

Everything has been obtained previously.

\section{More proofs and perspectives on co-Poisson}

The intertwining property was proven as a result on the adeles
and ideles, but obviously the proof can be written directly on
$\RR$. It will look like this, with $\varphi(y)$ an even Schwartz
function (and $\alpha(y)$ as above):

\begin{proof}
\From\  $\int_\RR \sum_{n\geq1}
|\varphi(ny)||\alpha(y)|\,dy<\infty$:
\begin{eqnarray*}
\int_\RR \sum_{n\geq1} \varphi(ny)\alpha(y)\,dy &=& \sum_{n\geq1}
\int_\RR \varphi(ny)\alpha(y)\,dy\\ =\sum_{n\geq1} \int_\RR
\varphi(y)\frac{\alpha(y/n)}{n}\,dy &=&\int_\RR
\varphi(y)\sum_{n\geq1} \frac{\alpha(y/n)}{n}\,dy\\
\end{eqnarray*}
On the other hand applying the usual Poisson summation formula:
\begin{eqnarray*}
&&\int_\RR \sum_{n\geq1} \varphi(ny)\alpha(y)\,dy\\ &=& \int_\RR
\left(\sum_{n\geq1} \frac{\cF(\varphi)(n/y)}{|y|} -
\frac{\varphi(0)}{2} +
\frac{\cF(\varphi)(0)}{2|y|}\right)\alpha(y)\,dy\\ &=&  \int_\RR
\Big(\sum_{n\geq1}
\cF(\varphi)(ny)\Big)\frac{\alpha(1/y)}{|y|}\,dy
-\varphi(0)\int_0^\infty \alpha(y)\,dy +
\cF(\varphi)(0)\int_0^\infty \frac{\alpha(y)}{|y|}\,dy \\ &=&
\int_\RR \cF(\varphi)(y)\sum_{n\geq1} \frac{\alpha(n/y)}{|y|}\,dy
-\varphi(0)\int_0^\infty \alpha(y)\,dy +
\cF(\varphi)(0)\int_0^\infty \frac{\alpha(y)}{|y|}\,dy \\
\end{eqnarray*}
The conclusion being:
$$\int_\RR \varphi(y)\sum_{n\geq1} \frac{\alpha(y/n)}{n}\,dy -
\cF(\varphi)(0)\int_0^\infty \frac{\alpha(1/y)}{y}\,dy = \int_\RR
\cF(\varphi)(y)\sum_{n\geq1} \frac{\alpha(n/y)}{|y|}\,dy
-\varphi(0)\int_0^\infty \alpha(y)\,dy$$
which, after exchanging $\varphi$ with $\cF(\varphi)$, is a
distribution theoretic formulation of the intertwining property:
$$\cF\left(\sum_{n\geq1} \frac{\alpha(y/n)}{n} - \int_0^\infty
\frac{\alpha(y)}{y}dy\right) = \sum_{n\geq1}
\frac{\alpha(n/y)}{|y|} - \int_0^\infty {\alpha(y)}dy$$
\end{proof}

It is useful to have a version of co-Poisson intertwining without
compactness nor smoothness conditions:

\begin{lem}
Let $g(u)$ be an even measurable function with
$$\int_0^\infty {|g(u)|\over u}du < \infty$$ The sum
$\sum_{n\geq1} {g(t/n)\over n}$ is Lebesgue almost-everywhere
absolutely convergent. It is a locally integrable function of
$t$. It is a tempered distribution.
\end{lem}

\begin{proof}
Let $\varphi(t)$ be an arbitrary even measurable function. One
has:
$$\int_0^\infty \sum_{n\geq1} |\varphi(nt)||g(t)|dt =
\int_0^\infty |\varphi(t)|\sum_{n\geq1} \left|{g(t/n)\over
n}\right|dt$$ If we take $\varphi(t)$ to be $1$ for $|t|\leq
\Lambda$, $0$ for $|t|>\Lambda$, we have $\sum_{n\geq1}
|\varphi(nt)| = [\Lambda/|t|] \leq \Lambda/|t|$.  From this:
$$\int_0^\Lambda \sum_{n\geq1} \left|{g(t/n)\over n}\right|dt =
O(\Lambda)$$
and this implies that $A(t) = \sum_{n\geq1}
{g(t/n)\over n}$ is almost everywhere absolutely convergent, that
it is locally integrable and also that $\int_0^u A(t)dt$ is
$O(u)$. So the continuous function $\int_0^u A(t)dt$ is a
tempered distribution.  Hence its distributional derivative
$A(u)$ is again a tempered distribution.
\end{proof}

\begin{theorem}\label{copoissongeneral}
Let $g(t)$ be an even measurable function with
$$\int_0^\infty {|g(t)|\over t}dt + \int_0^\infty |g(t)|dt <
\infty$$ Then the co-Poisson intertwining
$${\cal F}\left(\sum_{n\geq1} {g(t/n)\over n} - \int_0^\infty
{g(u)\over u}du\right) = \sum_{n\geq1} {g(n/t)\over t} -
\int_0^\infty {g(u)}du$$ holds as an identity of tempered
distributions.
\end{theorem}

\begin{proof}
The proof \ref{copoisson} given above applies identically. To get
it started one only has to state trivially for $\varphi(y)$ a
Schwartz function that $\sum_{n\geq1} |\varphi(ny)|$ is
$O(1/|y|)$.
\end{proof}

So the conclusion is that as soon as the two integrals are
absolutely convergent the full co-Poisson intertwining makes
sense and holds true. If one manages to get more information, the
meaning of the ${\cal F}$ will be accordingly improved. For
example if one side is in $L^2$ then the other side has to be too
and the equality holds with $\cF$ being the Fourier-Plancherel
(cosine) transform.

\begin{note}
In this familiar $\RR$ setting our first modification of the
Poisson formula was very cosmetic: the Poisson formula told us
the equality of two functions and we exchanged a term on the left
with a term on the right. This was to stay in a Hilbert space,
but it remained a statement about the equality of two functions
(and in the adelic setting, the equality of two functions with
values in the distributions on the  adeles). With the co-Poisson
if we were to similarly exchange the integral on the left with
the integral on the right, we would have to use Dirac
distributions, and the nature of the identity would change. So
the co-Poisson is more demanding than the (modified) Poisson.
Going from Poisson to co-Poisson can be done in many ways:
conjugation with $I$, or conjugation with $\cF$, or Hilbert
adjoint, or more striking still and at the same time imposed upon
us from adeles and ideles, the switch from viewing a certain
quantity as a function (Poisson) to viewing it as a distribution
(co-Poisson). The co-Poisson is a distribution whose values we
try to understand as $L^2$-functions, whereas the (modified)
Poisson is a function with values in distributions (bad for
Hilbert space).
\end{note}

We state again the important intertwining property as a theorem,
with an  alternative proof:

\begin{theorem}
Let $\alpha(y)$ be a smooth even function on $\RR$ with compact
support away from the origin. Let $P^\prime(\alpha)$ be its
co-Poisson summation:
$$P^\prime(\alpha)(y) = \sum_{n\geq1} \frac{\alpha(y/n)}{n} -
\int_0^\infty \frac{\alpha(y)}{y}\,dy$$
Then the additive Fourier Transform of $P^\prime(\alpha)$ is
$P^\prime(I(\alpha))$ with $I(\alpha)(y) = \alpha(1/y)/|y|$.
\end{theorem}

\begin{proof}
Let $P$ be the (modified) Poisson summation on even functions:
$$P(\alpha)(y) = \sum_{n\geq1} \alpha(ny) - \frac{\int_0^\infty
\alpha(y)\,dy}{|y|}$$
Obviously $P^\prime = I\cdot P\cdot I$. And we want to prove $\cF
P^\prime = P^\prime I$. Let us give a formal operator proof:
$$\cF P^\prime = \cF I P I = P \cF I I = P \cF = I  P = P^\prime
I$$
Apart from the usual Poisson summation formula $P\cF = I P$ the
crucial step was the commutativity of $\cF I$ and $P$. This
follows from the fact that both operators commute with the
multiplicative action of $\RR^\times$, so they are simultaneously
diagonalized by multiplicative characters, hence they have to
commute.

To elucidate this in a simple manner we extend our operators $I$,
$\cF$, $P$ and $P^\prime$ to a larger class of functions, a class
stable under all four operators. It is not difficult \cite{jnt}
to show that for each Schwartz function $\beta$ (in particular
for $\beta = I(\alpha)$) the Mellin Transform of $P(\beta)(y)$ is:
$$\int_0^\infty P(\beta)(y)y^{s-1}\,dy = \zeta(s)\int_0^\infty
\beta(y)y^{s-1}\,dy$$
 initially at least for $0<\Real(s)<1$. Let us consider the class
 of functions $k(1/2 + i\tau)$ on the critical line which
 decrease faster than any inverse power of $\tau$ when
 $|\tau|\to\infty$. On this class of functions we define $I$ as
 $k(s)\mapsto k(1-s)$, $P$ as $k(s)\mapsto \zeta(s)k(s)$,
 $P^\prime$ as $k(s)\mapsto \zeta(1-s)k(s)$, and $\cF I$ as
 $k(s)\mapsto \gamma_+(s) k(s)$ with $\gamma_+(s) = \pi^{-(s-
 {1/2})}{\Gamma(s/2)/\Gamma((1-s)/2)}$. The very crude bound (on
 $\Real(s) = 1/2$) $|\zeta(s)| = O(|s|)$ (for example, from
 $\zeta(s)/s = 1/(s-1) - \int_0^1 \{1/t\}t^{s-1}\,dt$) shows that
 it is a multiplier of this class (it is also a multiplier of the
 Schwartz class from the similar crude bounds on its derivatives
 one obtains from the just given formula). And $|\gamma_+(s)|=1$,
 so this works for it too (and also for the Schwartz class, see
 \cite{conductor}). The above formal operator proof is now not
 formal anymore (using, obviously, that the Mellin transform is
 one-to-one on our $\alpha$'s). The intertwining property for
 $P^\prime$ is equivalent to the intertwining property for $P$,
 because both are equivalent, but in different ways, to the
 functional equation for the zeta function.
\end{proof}

As was stated in the previous proof a function space which is
stable under all four operators $I$, $\cF$, $P$ and $P^\prime$ is
the space of inverse Mellin transforms of Schwartz functions on
the critical line: these are exactly the even functions on $\RR$
with the form $k(\log|y|)/\sqrt{|y|}$ where $k(a)$ is a Schwartz
function of $a\in\RR$. We pointed out the stability under Fourier
Transform in \cite{conductor}. We note that although $P$ and
$P^\prime$ make sense when applied to $k(\log|y|)/\sqrt{|y|}$ and
that they give a new function of this type, this can not always
be expressed as in their original definitions, for example
because the integrals involved have no reason to be convergent
(morally they correspond to evaluations away from the critical
line at $0$ and at $1$.)

We may also study $I$, $\cF$, $P$ and $P^\prime$ as operators on
$L^2$, but some minimal care has to be taken because $P$ and its
adjoint $P^\prime$ are not bounded. Nevertheless they are closed
operators and they commute with the Abelian von Neumann algebra
of bounded operators commuting with $\RR^\times$ (see
\cite{mrl}). This gives one more method to establish the
co-Poisson intertwining as a corollary to the Poisson
intertwining, as we may go from $P\cF = I P$ (\emph{modified}
Poisson) to $\cF P^\prime = P^\prime I$ (co-Poisson) simply by
taking Hilbert adjoints: $P^\prime = P^*$ (and here $\cF$ is the
cosine transform, so $\cF^* = \cF$):
$$P\cF = I P\ \Rightarrow\ (P\cF)^* = (I P)^*\ \Rightarrow\ \cF
P^\prime = P^\prime I$$
It is important to be aware that for this to be correct  it is
crucial that we are using $P$ to denote, not the \emph{original}
Poisson sums, but the \emph{M\"untz-modified} Poisson sums.

Another perspective on co-Poisson comes from a re-examination of
the use of the multiplicative version (now called Mellin) of the
Fourier Transform in Riemann's paper. To establish the functional
equation, Riemann uses the \emph{left} Mellin transform
$\int_0^\infty f(u) u^{s-1}du$, but when he relates the zeros to
the primes with an explicit formula he uses the \emph{right}
Mellin transform $\int_0^\infty f(u) u^{-s}du$. After all the
zeta-function itself is most simply expressed as the right Mellin
transform of the sum of the Dirac at the positive integers. We
said earlier that the M\"untz-modified sums corresponded under
Mellin Transform to the zeta function $\zeta(s)$, but this is
using the \emph{left} convention. If, rather, we use the
\emph{right} convention we are bound to associate to $\zeta(s)$
the \emph{co-Poisson sums}! If we now ask what the functional
equation tells us, then the answer is: in particular the
co-Poisson intertwining\dots It should be clear from this
discussion that the co-Poisson intertwining is a formula of the
nineteenth century which was discovered at the close of the
twentieth century.

Immediately after being communicated the co-Poisson formula, Luis
B\'aez-Duarte replied that an application of Euler-McLaurin
summation establishes the co-Poisson intertwining in a more
elementary manner, inasmuch as his method uses neither
distributions nor Mellin transforms, and does not appeal directly
to either the Poisson summation formula, nor to the functional
equation of the Riemann zeta function. Here is the proof emerging
from this discussion:

\begin{theorem}[proven jointly with Luis B\'aez-Duarte]
Let $g(t)$ be an (even) function of class $C^2$ which has compact
support away from the origin. Then both $\sum_{n\geq1}
{g(t/n)\over n} - \int_0^\infty {g(u)\over u}du$ and
$\sum_{n\geq1} {g(n/t)\over t} - \int_0^\infty {g(u)}du$ are
continuous $L^1$-functions and the co-Poisson intertwining formula
$${\cal F}\left(\sum_{n\geq1} {g(t/n)\over n} - \int_0^\infty
{g(u)\over u}du\right) = \sum_{n\geq1} {g(n/t)\over t} -
\int_0^\infty {g(u)}du$$ holds as a pointwise equality and may be
established as a corollary to the Euler-McLaurin summation
formulae.
\end{theorem}

\begin{proof}
Let $B(u)$ be the periodic funtion which on $[0,1)$ has values
${u^2\over2} - {u\over2} + {1\over12}$. It is a continuous even
function, which as is well-known is also expressed as:

$$B(u) = \sum_{n\geq1} {\cos(2\pi n u)\over2\pi^2 n^2}$$

Let $f(t)$ be a $C^2$ function with compact support away from
$0$. We have:

$$\sum_{n\geq1} f(n) - \int_0^\infty f(u) du = - \int_0^\infty
B(u)\left({d\over du}\right)^2 f(u)\,du$$

If we apply this to the function $u\to f(u/w)/w$ for $w>0$, we
obtain:

$$\sum_{n\geq1} {f(n/w)\over w} - \int_0^\infty f(u) du = -
\int_0^\infty B(u)\left({d\over du}\right)^2 {f(u/w)\over w}\,du
= -\frac1{w^2}\int_0^\infty B(wv)\left({d\over dv}\right)^2
{f(v)}\,dv$$

The left-hand side is locally a finite sum, hence of class $C^2$
(and when $|w|$ is sufficiently small it  reduces to the constant
$- \int_0^\infty f(u) du$) and the formula above shows that it is
$O(1/w^2)$ when $w\to\infty$. If we only assume $f(u)$ to be
bounded but still with compact support away from the origin then
obviously $\sum_{n\geq1} {f(n/w)\over w}$ is at any rate bounded
(the number of non-vanishing terms being $O(w)$). We now apply
the above formula to $g(t) = f(1/t)/t$:

\begin{eqnarray*}
\sum_{n\geq1} {g(t/n)\over n} - \int_0^\infty {g(1/u)\over u}du
&=& -\int_0^\infty B(u)\left({d\over du}\right)^2 {g(t/u)\over
u}\,du \\ = -\sum_{n\geq1} \int_0^\infty {\cos(2\pi n
u)\over2\pi^2 n^2}\left({d\over du}\right)^2 {g(t/u)\over u}\,du
&=& -\sum_{n\geq1} \int_0^\infty {\cos(2\pi tw)\over2\pi^2
t^2}\left({d\over dw}\right)^2 {g(n/w)\over w}\,dw \\
\end{eqnarray*}

At this stage we expand the second derivative and using that
$\sum_{n\geq1} {k(n/w)\over w}$ is bounded with $k(t) = |g(t)|,
|tg^\prime(t)|, |t^2 g^{\prime\prime}(t)|$ we see that dominated
convergence applies. So:

\begin{eqnarray*}
&=&-\int_0^\infty {\cos(2\pi tw)\over2\pi^2 t^2}\left({d\over
dw}\right)^2 (\sum_{n\geq1} {g(n/w)\over w})\,dw \\
&=&-\int_0^\infty {\cos(2\pi tw)\over2\pi^2 t^2}\left({d\over
dw}\right)^2 (\sum_{n\geq1} {g(n/w)\over w} - \int_0^\infty
g(\alpha)d\alpha)\,dw \\ &=&+\int_0^\infty {2\cos(2\pi
tw)}(\sum_{n\geq1} {g(n/w)\over w} - \int_0^\infty
g(\alpha)d\alpha)\,dw \\
\end{eqnarray*}

The first integration by parts at the end is justified by
$\lim_{w\to\infty} {d\over dw}\sum_{n\geq1} {g(n/w)\over w} = 0$,
using that $\sum_{n\geq1} {k(n/w)\over w}$ is bounded with $k(t)
= |g(t)|, |tg^\prime(t)|$. So far we have $t\neq 0$, but the
integral being dominated we may let $t=0$ in the final
formula. If we now replace $g(t)$ with $g(1/t)/t$ we get the
co-Poisson intertwining as a pointwise equality.
\end{proof}

The most direct attack when first confronted with the co-Poisson
formula is presumably this:  $\cF(\sum g(t/n)/n) = \sum
\cF(g)(n\xi) = \sum g(m/\xi)/|\xi|$, where we first intervert and
then apply Poisson. This has a number of pitfalls (for whose
unraveling the language of distributions is very useful), but it
is possible to make it work to prove something correct. This
proof is a joint-work with Bernard Candelpergher, and as in the
previous approach it uses only elementary tools and especially
the (simplest cases of) Euler-McLaurin summation.

\begin{theorem}[proven jointly with Bernard Candelpergher]
Let $g(t)$ be an even function of class $C^2$ which has compact
support away from the origin. Then both $\sum_{n\geq1}
{g(t/n)\over n} - \int_0^\infty {g(u)\over u}du$ and
$\sum_{n\geq1} {g(n/t)\over t} - \int_0^\infty {g(u)}du$ are
continuous $L^1$-functions and the co-Poisson intertwining formula
$${\cal F}\left(\sum_{n\geq1} {g(t/n)\over n} - \int_0^\infty
{g(u)\over u}du\right) = \sum_{n\geq1} {g(n/t)\over t} -
\int_0^\infty {g(u)}du$$ holds as a pointwise equality and may be
established as a corollary to the Euler-McLaurin summation
formulae.
\end{theorem}

\begin{proof}
We work on $(0,\infty)$. We have ($t>0$):
$$\sum_{n\geq1} \frac{g(t/n)}n - \int_0^\infty \frac{g(1/u)}u\,du
= \sum_{n\geq1} A_n(t)$$

with:

$$A_n(t) =
\frac12\left(\frac{g(t/n)}n+\frac{g(t/(n-1))}{n-1}\right) -
\int_{n-1}^n  \frac{g(t/u)}u\,du$$

where the term with $n-1$ is dropped for $n=1$. With
$k(t)=g(1/t)/t$ one has:

$$A_n(t) = \frac12\left(\frac{k(n/t)}t+\frac{k((n-1)/t)}t\right)
- \int_{n-1}^n  \frac{k(u/t)}t\,du = -  \int_{n-1}^n
C_2(u){\frac{d^2}{du^2}}\frac{k(u/t)}t\,du$$

where $C_2(u) = (\{u\}^2 - \{u\})/2$. From this:

$$\sum_{n\geq1} |A_n(t)| \leq \frac1{t^2}\int_0^\infty
|C_2(u)||k^{\prime\prime}(\frac ut)|\frac{du}t = O(\frac1{t^2})$$

\let\wt=\widetilde

On the other hand for $t$ small enough one has $A_n(t)\equiv 0$
for $n\geq2$. So the sum $\sum_{n\geq1} A_n(t)$ is absolutely
convergent to an $L^1$-function and also we can compute its
Fourier transform termwise. This gives, with $\wt g = \cF(g)$:

$$\cF(\sum_{n\geq1} A_n(t))(\xi) = \frac{-\wt g(0)}2 +
\sum_{n\geq1} \left( \frac{\wt g(n\xi)+\wt g((n-1)\xi)}2  -
\int_{n-1}^n  \wt g(u\xi)\,du\right)$$

The appearance of $-\wt g(0)/2$ is from the fact that this time
the $n-1$ term with $n=1$ is to be counted in the sum, so we have
to compensate for this.  The formula holds for all $\xi$
($\xi\geq0$), in particular for $\xi = 0$ it reads:

$$\int_\RR \left(\sum_{n\geq1} \frac{g(t/n)}n - \int_0^\infty
\frac{g(1/u)}u\,du\right)\,dt = - \int_0^\infty g(v)\,dv$$

For $\xi>0$ we are one step away from co-Poisson, it only remains
to apply Poisson to our sum (note that $\int_0^\infty \wt
g(u\xi)\,du = 0$). But we can also base this on
Euler-McLaurin. Indeed with $B_1(v) = \{v\} - \frac12$:

$$\frac{\wt g(n\xi)+\wt g((n-1)\xi)}2  - \int_{n-1}^n  \wt
g(u\xi)\,du = \int_{n-1}^n B_1(v)\frac{d}{dv}\wt g(v\xi)\,dv$$

Hence

$$\cF(\sum_{n\geq1} A_n(t))(\xi) =  - \int_0^\infty g(v)\,dv +
\int_0^\infty B_1(v)\frac{d}{dv}\wt g(v\xi)\,dv$$

Now, from the fact that $g$ is $C^2$ with compact support, its
Fourier transform and all derivatives of it are
$O(1/|\xi|^2)$. We may thus use the boundedly convergent
expression:

$$B_1(v) = \sum_{m\geq1} \frac{-\sin(2\pi mv)}{m\pi}$$

and then

\begin{eqnarray*}
\int_0^\infty B_1(v)\frac{d}{dv}\wt g(v\xi)\,dv &=&
\lim_{M\to\infty}\sum_{1\leq m\leq M}
\int_0^\infty\frac{-\sin(2\pi mv)}{m\pi} \frac{d}{dv}\wt
g(v\xi)\,dv \\ &=& \lim_{M\to\infty}\sum_{1\leq m\leq M}
2\int_0^\infty \cos(2\pi mv) \wt g(v\xi)\,dv \\ &=&
\lim_{M\to\infty}\sum_{1\leq m\leq M}  \frac{g(m/\xi)}\xi\\
\end{eqnarray*}

where Fourier inversion was used ($\xi>0$). This completes the
proof of co-Poisson. It is also interesting to prove in another
manner the special formula ($\xi = 0$):

$$\int_0^\infty \left(\sum_{n\geq1} \frac{g(t/n)}n -
\int_0^\infty \frac{g(1/u)}u\,du\right)\,dt = - \frac12
\int_0^\infty g(v)\,dv$$

This may be done as follows: first,

$$\int_0^\Lambda \sum_{n\geq1} \frac{g(t/n)}n dt = \sum_{n\geq1}
\int_0^{\Lambda/n} g(t)dt = \int_0^\infty \left[\frac\Lambda
t\right]g(t)dt$$

so we are looking at

$$\lim_{\Lambda\to\infty} -\int_0^\infty\left\{\frac\Lambda
t\right\}g(t)dt = -\lim_{\Lambda\to\infty} \int_0^\infty
\{\Lambda v\}h(v)dv$$

with the $L^1$-function $h(v) = g(1/v)/v^2$. It is obvious that

$$0\leq A\leq B\quad\Rightarrow\quad \lim_{\Lambda\to\infty}
\int_0^\infty \{\Lambda v\}\Un_{A\leq v\leq B}(v)dv =
\frac{B-A}2$$

so using the density argument from the usual proof of the
Riemann-Lebesgue lemma one deduces that

$$\lim_{\Lambda\to\infty} \int_0^\infty \{\Lambda v\}h(v)dv =
\frac12\int_0^\infty h(v)dv$$

for all $L^1$-functions. One last remark is that at the level of
(left) Mellin, co-Poisson is like multiplication by $\zeta(1-s)$,
so the special formula is just another instance of $\zeta(0) =
-\frac12$.
\end{proof}

Here is one last approach to co-Poisson (extracted from a
manuscript in preparation \cite{copoisson}). The Poisson
summation identity is
$$\sum_{n\in\ZZ} \cF(f)(n) = \sum_{m\in\ZZ} f(m)$$ It applies in
particular to Schwartz functions, and may be written as:
$$\cF\left(\sum_{n\in\ZZ} \delta_n(x)\right) = \sum_{m\in\ZZ}
\delta_m(y)$$ We take this identity of tempered distributions
seriously, on its own, and do not completely identify it with the
Poisson summation identity. Let $t\neq0$ and let us replace $x$
by $tx$. One has $\delta_n(tx) = \delta_{n/t}(x)/|t|$. So
$$\cF\left(\sum_{n\in\ZZ} \frac{\delta_{n/t}(x)}{|t|}\right) =
\frac1{|t|}\sum_{m\in\ZZ} \delta_m(\frac{y}{t}) = \sum_{m\in\ZZ}
\delta_{tm}(y)$$
We average these tempered distributions with an integrable weight
$g(t)$ (which is compactly supported away from the origin) to
obtain an identity of tempered distributions:
$$\cF\left(\sum_{n\in\ZZ} \int
g(t)\frac{\delta_{n/t}(x)}{|t|}\,dt\right) = \sum_{m\in\ZZ} \int
g(t)\delta_{tm}(y)\,dt$$
But for $n\neq0$ (resp. $m\neq0$) and as distributions in $x$
(resp. $y$):
$$\int g(t)\frac{\delta_{n/t}(x)}{|t|}\,dt = \int
g(\frac{n}\alpha)\delta_\alpha(x)\,\frac{d\alpha}{|\alpha|}
=\frac{g(\frac{n}x)}{|x|}$$
$$\int g(t)\delta_{tm}(y)\,dt = \int
g(\frac\beta{m})\delta_\beta(y)\,\frac{d\beta}{|m|}
=\frac{g(\frac{y}{m})}{|m|}$$
whereas $\int g(t)\frac{\delta_{0}(x)}{|t|}\,dt = (\int
g(t)\frac{dt}{|t|})\delta_{0}(x)$ and $\int g(t)\delta_{0}(y)\,dt
= (\int g(t)dt)\delta_{0}(y)$. If we exchange sides for the
contributions of $n=0$ and $m=0$ we end up with the co-Poisson
intertwining as an identity of tempered distributions. This
method of multiplicative convolution applies to situations where
the discreteness of the support of the original distributions
applies only at the origin. A general discussion is planned in
\cite{copoisson}.

There is reason to believe that the problem of understanding the
spaces of $L^2$-functions which are vanishing together with their
Fourier Transform in a neighborhood of the origin, is important
simultaneously for Analysis and Arithmetic. It is a remarkable
ancient discovery of de~Branges \cite{Bra64} \cite{Bra} that
these spaces have the rich structure which he developed generally
in his theory of Hilbert Spaces of entire functions: they are
among the ``Sonine spaces''. The co-Poisson summations have the
(extended) property of being constant together with their Fourier
transform, in some neighborhood of the origin, and we will show
later (see also \cite{twosys}) that the zeros of the Riemann zeta
function are the obstructions for (the square-integrable among)
the co-Poisson sums to fill up the full spaces of such
square-integrable functions.

\section{Impact of the B\'aez-Duarte, Balazard, Landreau
  and Saias theorem on the so-called Hilbert-P\'olya idea}
  
How could it be important to replace $\zeta(s)$ with
$\zeta(1-s)$\thinspace? Clearly only if we leave the critical
line and start paying attention to the difference between the
right half-plane $\Real(s)>1/2$ and the left half-plane
$\Real(s)<1/2$. Equivalently if we switch from the full group of
contractions-dilations $C_\lambda$, $0<\lambda<\infty$, which
acts as $\phi(x)\mapsto \phi(x/\lambda)/\sqrt{|\lambda|}$ on even
functions on $\RR$, or as $Z(s)\mapsto \lambda^{s-1/2}Z(s)$ on
their Mellin Transforms, to its sub-semi-group of contractions
($0<\lambda\leq1$). The contractions act as isometries on the
Hardy space $\HH^2(\Real(s)>1/2)$ or equivalently on its inverse
Mellin transform the space $L^2((0,1), dt)$.

It is a theme contemporaneous to Tate's Thesis and Weil's first
paper on the explicit formula that it is possible to formulate
the Riemann Hypothesis in such a semi-group set-up: this is due
to Nyman \cite{Nym} and Beurling \cite{Beu55} and builds on the
Beurling \cite{Beu49} (and later for the half-plane) Lax
\cite{Lax59} theory of invariant subspaces of the Hardy
spaces. The criterion of Nyman reads as follows: the linear
combinations of functions $t\mapsto \{1/t\} - \{a/t\}/a$, for
$0<a<1$, are dense in $L^2((0,1), dt)$ if and only if the Riemann
Hypothesis holds. It is easily seen that the smallest closed
subspace containing these functions is stable under contractions,
so to test the closure property it is only necessary to decide
whether the constant function $\Un$ on $(0,1)$ may be
approximated. The connection with the zeta function is
established through ($0<\Real(s)<1$):
$$\int_0^\infty \{\frac{1}{t}\}t^{s-1}dt = - \frac{\zeta(s)}{s}$$
This gives for our functions, with $0<a<1$ and $0<\Real(s)$:
$$\int_0^1 (\{\frac{1}{t}\} - \frac{1}{a}\{\frac{a}{t}\})
t^{s-1}dt = (a^{s-1} - 1)\frac{\zeta(s)}{s}$$
The question is whether the invariant (under contractions)
subspace of $\HH^2(\Real(s)>1/2)$ of linear combinations of these
Mellin transforms is dense or not. Each zero $\rho$ of $\zeta(s)$
in $\Real(s)>1/2$ is an obvious obstruction as (the complex
conjugate of) $t^{\rho -1}$ belongs to $L^2((0,1), dt)$. The
Beurling-Lax theory describing the structure of invariant
subspaces allows the conclusion in that case that there are no
other obstructions (we showed \cite{ams} as an addendum that the
norm of the orthogonal projection of $\Un$ to the Nyman space is
$\prod_{\Real(\rho)>1/2} |(1-\rho)/\rho|$, where the zeros are
counted with their multiplicities).  Recently
(\cite{luisnatural}) Luis B\'aez-Duarte has shown that the appeal
to the Beurling-Lax invariant subspace theory could be completely
avoided, and furthermore he has put the Nyman-Beurling criterion
in the stronger form where one applies  to the fractional part
function only integer-ratio contractions.

One could think from our description of the original proof of the
Nyman-Beurling criterion that the zeros \emph{on the critical
line} are out of its scope, as they don't seem to play any
r\^ole. So it has been a very novel thing when B\'aez-Duarte,
Balazard, Landreau and Saias asked the right question and
provided a far from obvious answer \cite{Bal3}. First a minor
variation is to replace the Nyman criterion with the question
whether the function $\Un_{0<t<1}$ can be approximated in
$L^2(0,\infty)$ with linear combinations of the contractions of
$\{1/t\}$. Let $D(\lambda)$ be the Hilbert space distance between
$\Un_{0<t<1}$ and linear combinations of contractions
$C_\theta(\{1/t\})$, $\lambda\leq \theta\leq1$. The Riemann
Hypothesis holds if and only if $\lim D(\lambda) = 0$.

{\bf Theorem of B\'aez-Duarte, Balazard, Landreau and Saias
\cite{Bal3}: } \emph{ One has the lower bound:}
$$\liminf\; |\log(\lambda)|\;D(\lambda)^2 \geq \sum_\rho
\frac{1}{|\rho|^2}$$
\emph{where the sum is over all non-trivial zeros of the zeta
function, counted only once independently of their multiplicity.}

The authors of \cite{Bal3} conjecture that  equality holds (also
with $\lim$ in place of $\liminf$) when one counts the zeros with
their multiplicities: our next result shows that their conjecture
not only implies the Riemann Hypothesis but it also implies the
simplicity of all the zeros:

\begin{theorem}[\cite{aim}]\label{jfapprox}
The following lower bound holds:
$$\liminf\; |\log(\lambda)|\;D(\lambda)^2 \geq \sum_\rho
\frac{m_\rho^2}{|\rho|^2}$$
\end{theorem}

Our proof relies on the link we have established between the
study originated by B\'aez-Duarte, Balazard, Landreau and Saias
of the distance function $D(\lambda)$ and the so-called
Hilbert-P\'olya idea. This idea will be taken here in the
somewhat vague acceptation that the zeros of $L$-functions may
have a natural interpretation as Hilbert space vectors,
eigenvectors for a certain self-adjoint operator. If we had such
vectors in $L^2((0,\infty), dt)$, perpendicular to $\{1/t\}$ and
to its contractions $C_\theta(\{1/t\})$, $\lambda\leq
\theta\leq1$ then we would be in position to obtain a lower bound
for $D(\lambda)$ from the orthogonal projection of $\Un_{0<t<1}$
to the space spanned by the vectors. This lower bound would be
presumably easily expressed as a sum indexed by the zeros from
the fact that eigenspaces of a self-adjoint operator are mutually
perpendicular. The first candidates are $t\mapsto t^{-\rho}$:
they satisfy formally the perpendicularity condition to $\{1/t\}$
and its contractions, but they do not belong to
$L^2$. Nevertheless we could be in a position to approximately
implement the idea if we used instead the square integrable
vectors $t\mapsto t^{-(\rho-\epsilon)}\Un_{0<t<1}$, $\epsilon>0$,
$\Real(\rho)=1/2$. The authors of \cite{Bal3} followed more or
less this strategy, but as they did not benefit from exact
perpendicularity, they had to provide not so easily obtained
estimates. It appears that the $\epsilon>0$ does not seem to
allow to take easily into account the multiplicities of the
zeros. For their technical estimates the authors of \cite{Bal3}
used to great advantage a certain scale invariant operator $U$,
which had been introduced by B\'aez-Duarte in an earlier paper
\cite{Bae} discussing the Nyman-Beurling problem.

How is it possible that the use of this B\'aez-Duarte operator
$U$ (whose definition only relies on some ideas of harmonic
analysis, and some useful integral formulae, with at first sight
no arithmetic involved) allows B\'aez-Duarte, Balazard, Landreau
and Saias to make progress on the Nyman-Beurling
criterion\thinspace? We related the mechanism underlying the
insightful B\'aez-Duarte construction \cite{Bae} of the operator
$U$ to a construction quite natural in scattering theory and it
appeared then that it was possible to use it (or a variant $V$)
to construct true Hilbert space vectors $Y^\lambda_{s,k}$ indexed
by $s$ on the critical line, and $k\in\NN$, having the property
of expressing the values of the Riemann zeta function and its
derivatives on the critical line as Hilbert space scalar products:
$$(A, Y^\lambda_{s,k}) = {(-\frac{d}{ds})}^k\;
\frac{s-1}{s}\frac{\zeta(s)}{s}$$
The additional factors are such that
$\frac{s-1}{s}\frac{\zeta(s)}{s}$ belongs to the Hardy space of
the right half-plane and $A(t)$ is its inverse Mellin transform,
an element of $L^2((0,1),dt)$. What is more, we can replace $A$
with its contractions $C_\theta(A)$ as long as $\lambda\leq
\theta \leq 1$:
$$\lambda\leq \theta \leq 1 \Rightarrow (C_\theta(A),
Y^\lambda_{s,k}) = {(-\frac{d}{ds})}^k\;\theta^{s-1/2}\;
\frac{s-1}{s}\frac{\zeta(s)}{s}$$

So the only ones among the $Y^\lambda_{s,k}$ which are (exactly,
not approximately) perpendicular to $A$ and its contractions up
to $\lambda$ are the $Y^\lambda_{\rho,k}$, $\zeta(\rho) = 0$,
$k<m_\rho$ ($\lambda<1$). This connects the Nyman-Beurling
criterion with the so-called Hilbert-P\'olya idea.  From  the
explicit integral formulae of B\'aez-Duarte for his operator $U$
one can write formulae for the vectors $Y^\lambda_{s,k}$ from
which the following asymptotic behavior emerges: \let\la=\lambda
\begin{eqnarray*}
\lim_{\la\to0}\ |\log(\lambda)|^{-1-k-l}\cdot(Y^\la_{s_1,k},
Y^\la_{s_2,l}) &=& 0\qquad(s_1\neq s_2)\\ \lim_{\la\to0}\
|\log(\lambda)|^{-1-k-l}\cdot(Y^\la_{s,k}, Y^\la_{s,l}) &=&
{1\over k+l+1}
\end{eqnarray*}
In particular the rescaled vectors $X^\lambda_{s,0} =
Y^\lambda_{s,0}/\sqrt{|\log(\lambda)|}$ become orthonormal in the
limit when $\lambda\to0$. The limit can not work inside
$L^2((0,\infty))$ because this is a separable space\thinspace! In
fact one shows without difficulty that the vectors
$X^\lambda_{s,k}$ weakly converge to $0$ as  $\lambda\to0$. The
theorem \ref{jfapprox} is easily deduced from the above estimates.

We have not explained yet what is $U$ and how the
$Y^\lambda_{s,k}$ are constructed with the help of it. The
B\'aez-Duarte operator is the unique scale invariant operator
with sends $\{1/t\}$ to its image under the \emph{time reversal}
$J:\varphi(t)\mapsto \overline{\varphi(1/t)/t}$, which is here
the function $\{t\}/t$. Equivalently at the level of the Mellin
Transforms, $U$ acts as a multiplicator with multiplier on the
critical line
$$U(s) = \frac{\ \overline{\zeta(s)/s}\ }{\zeta(s)/s} =
\frac{\zeta(1-s)}{\zeta(s)}\; \frac{s}{1-s}$$

So this construction is extremely general: as soon as the
function $A\in L^2((0,\infty),dt)$ has an almost everywhere
non-vanishing Mellin transform $Z(s)$ on the critical line (which
by a theorem of Wiener is  equivalent to the fact that the
multiplicative translates $C_\theta(A)$, $0<\theta<\infty$, span
$L^2$) then we may associate to it the scale invariant operator
$V$ with acts as $\overline{Z(s)}/Z(s)$, and sends $A$ to its
time reversal $J(A)$. We see that $V$ is necessarily unitary. Let
us in particular suppose that $A$ is in $L^2((0,1),dt)$: then
$J(A)$ has support in $[1,\infty)$ and the images under $V$ of
the contractions $C_\theta(A)$, $\lambda\leq\theta\leq1$, being
contractions of $J(A)$, will be in $L^2((\lambda,\infty),dt)$. In
the case at hand we have:

$$Z(s) = \frac{s-1}{s}\frac{\zeta(s)}{s}\qquad V(s) =
\frac{\zeta(1-s)}{\zeta(s)}{\left(\frac{s}{1-s}\right)}^3$$

The B\'aez-Duarte operator $U$ and its cousin $V$ depend on
$\zeta(s)$ only  through its functional equation, which means
that they are associated with the (even) Fourier transform
$\cF_+$ (the cosine transform). We can use (almost) the same
operators for Dirichlet $L$-series with an even character (with
due attention paid to the conductors $q>1$), and there are other
operators we would use for odd characters, associated with the
sine transform $\cF_-$.

The vectors $Y^\lambda_{s,k}$ are obtained as follows: we start
from $|\log(t)|^k\;t^{-(s-\epsilon)}\Un_{0<t<1}$, apply $V$,
restrict to $[\lambda, +\infty)$, take the limit which now exists
in $L^2$ as $\epsilon\to0$, and apply $V^{-1}$. What happens is
that $V(t^{-1/2 - i\tau}\Un_{0<t<1})$ does not belong to $L^2$
but this is entirely due to its singularity at $0$, which, it
turns out, is $V(1/2+i\tau)t^{-1/2 - i\tau}$. We would not expect
it to be possible that a localized singularity would remain
localized after the  action of $\cF_+$ but the point is that the
operator with multiplier $\zeta(1-s)/\zeta(s)$ is the composite
$\cF_+\cdot I$ with $I(\phi)(x) = \phi(1/x)/|x|$. So the
singularity is first sent to infinity and $\cF_+$ puts it back at
the origin.

\begin{note}
Let us denote by $L$ the scale invariant unitary operator with
spectral function $L(s) = s/(s-1)$. One has $L = 1 - M$ where $M$
is the Hardy averaging $\phi(t)\mapsto(\int_0^t
\phi(u)\,du)/t$. The operator $V$ is $(-L)^3\cF_+ I = \cF_+ I
(-L)^3$. One has $ILI = L^{-1}$ and $LL^* = 1$. The operator $L$
is ``real'', meaning that it commutes  with the anti-unitary
complex conjugation $\phi(t)\mapsto \overline{\phi(t)}$. One has
$IV = I\cF_+ I (-L)^3$. We will write $\cG_+ = I\cF_+ I$, so that
$IV = \cG_+ (-L)^3 = (-L)^{-3}\cG_+$. The operator $\cG_+$ is
unitary and satisfies $\cG_+^2 =  1$. The operator $\cG_+$ is
real. The $U$ operator of B\'aez-Duarte is $\cF_+ I(-L)$, so $V =
UL^2 = L^2 U$. The operators $I$, $\cF_+$, $\cG_+$, $U$ and $V$
are real.
\end{note}

We now proceed with a more detailed study of the vectors
$Y^\lambda_{s,k}$, and of their use  to express values of Mellin
transforms and their derivatives on and off the critical line as
Hilbert space scalar products. We established
$$(B, Y^\lambda_{s,k}) = {(-\frac{d}{ds})}^k\;\widehat{B}(s)$$
for the Hardy functions $C_\theta(A)$,
$\lambda\leq\theta\leq1$. Their Mellin Transforms $\int_0^\infty
C_\theta(A)(t)t^{s-1}\,dt= \theta^{s-1/2}\;
\frac{s-1}{s}\frac{\zeta(s)}{s}$ are analytic in the entire
complex plane except for a double pole at $s=0$. The vectors
$Y^\lambda_{s,k}$ are the analytic continuation to $\Real(s)=1/2$
of vectors with the same definition  $Y^\lambda_{w,k}$,
$\Real(w)<1/2$. The above equation has thus its right hand side
analytic in $s$ but its left-hand side seemingly anti-analytic,
as our scalar product is linear in its first factor and conjugate
linear in its second factor. So we will use rather the
\emph{euclidean bilinear form} $[B,C] = \int_0^\infty
B(t)C(t)\,dt$. The spaces we consider are stable under complex
conjugation $B(t)\mapsto \overline{B(t)}$, and the operators we
use are real, so  statements of perpendicularity may equivalently
be stated using either $[\cdot,\cdot]$ or $(\cdot,\cdot)$. The
identity can then be restated for all  finite linear combinations
$B$ of our $C_\theta(A)$'s, $\lambda\leq\theta\leq1$, as
$$[B,Y^\lambda_{w,k}] = {(+\frac{d}{dw})}^k \int_0^1
t^{-w}B(t)\,dt = {(\frac{d}{dw})}^k
\left(\widehat{B}(1-w)\right)$$
for $\Real(w)\leq 1/2$. If we look at the proof of the main
theorem in \cite{aim} we see that the only thing that matters
about $B$ is that it should be supported in $[0,1]$ and that
$V(B)$ should be supported in $[\lambda,\infty)$, equivalently
that $(IV)(B)(t)$ has support in $[0,\Lambda]$, $\Lambda =
1/\lambda$. Let us note the following:

\begin{theorem}
The real unitary operator $IV$ satisfies $(IV)^2=1$.
\end{theorem}

\begin{note}
In particular  $IV$ is what B\'aez-Duarte calls ``a skew-root''
\cite{Bae}.
\end{note}

\begin{proof}
This is clear from the spectral representation
$$V(s) =
\frac{\zeta(1-s)}{\zeta(s)}{\left(\frac{s}{1-s}\right)}^3$$
which shows that $IVI = V^*$.
\end{proof}

We let $\Lambda = 1/\lambda$ and $\cG_\Lambda = IVC_\Lambda =
C_\lambda IV$. We note that $(\cG_\Lambda)^2 = 1$. We also note
$V\cG_\Lambda = C_\lambda VIV = C_\lambda I$. Let $M_\Lambda =
\HH^2\; \cap\; \cG_\Lambda(\HH^2)$, where we use the notation
$\HH^2 = L^2((0,1),dt)$. Obviously $\cG_\Lambda(M_\Lambda) =
M_\Lambda$. The function $A$ as well as its contractions
$C_\theta(A)$, $\lambda\leq\theta\leq1$ belong to
$M_\Lambda$. Indeed $\cG_\Lambda(A) = C_\lambda(A)$. We note that
the Mellin transform $\int_0^1 B(t) t^{w-1}\,dt$ of $B\in
M_\Lambda$ is analytic at least in $\Real(w)>1/2$. Let
$Q_\lambda$ be the orthogonal projection to $L^2(\lambda,\infty)$.

\begin{theorem}
The vectors $Y^\lambda_{w,k}$, originally defined for
$\Real(w)<1/2$ as
$$V^{-1}Q_\lambda V(|\log(t)|^k t^{-w}\Un_{0<t<1})$$
have (inside $L^2$) an analytic continuation in $w$  to the
entire complex plane $\CC$ except at $w=1$. The Mellin Transform
of $B\in M_\Lambda$ has an analytic continuation to
$\CC\setminus\{0\}$, with at most a pole of order $2$ at
$w=0$. One has for $w\neq1$ and $k\in \NN$:
$$[B,Y^\lambda_{w,k}] = {(\frac{d}{dw})}^k
\left(\widehat{B}(1-w)\right)$$
The following functional equation holds:
$$\widehat{\cG_\Lambda(B)}(w) =
\lambda^{w-1/2}V(1-w)\widehat{B}(1-w)$$
One has
$$\forall B\in M_\Lambda\quad \widehat{B}(-2) = \widehat{B}(-4) =
\cdots  = 0$$
\end{theorem}

\begin{proof}
We leave the details of the case $k>0$ to the reader. Let first
$\Real(w)>1/2$. We have for $B\in M_\Lambda$:
$$\widehat{B}(w) = \int_0^1 B(t)t^{w-1}\,dt = [B,
t^{w-1}\Un_{0<t<1}] =  [V(B),V(t^{w-1}\Un_{0<t<1})]$$
Writing $B = \cG_\Lambda(C)$, with $C\in\HH^2$, we get $V(B) =
C_\lambda I(C)$. So $V(B)$ has its support in $[\lambda,\infty)$
and:
$$\widehat{B}(w) = \int_\lambda^\infty
V(B)(u)V(t^{w-1}\Un_{0<t<1})(u)\,du$$
We will show that $V(t^{w-1}\Un_{0<t<1})(u)$ is analytic, for
fixed $u$, in $w\in\CC\setminus\{0\}$ and that it is
$O((1+|\log(u)|)/u)$ on $[\lambda,\infty)$, uniformly when $w$ is
in a compact subset of $\CC\setminus\{0\}$ (this is one logarithm
better than the estimate in \cite{aim} for $\Real(1-w)<1$). We
will thus have obtained the analytic continuation of the vectors
$Y^\lambda_{1-w,0}$ from $\Real(w)>1/2$ and at the same time the
analytic continuation of $\widehat{B}(w)$ as well as the formula:
$$w\neq0\ \Rightarrow\ [B,Y^\lambda_{1-w,0}] = \widehat{B}(w)$$

So the problem is to study the analytic continuation of
$V(t^{-z}\Un_{0<t<1})(u)$ from $\Real(z)<1/2$. If we followed the
method of \cite{aim}, we would write $V=(1-M)^2 U$, compute some
explicit formula for $U(t^{-z}\Un_{0<t<1})(u)$ and work with
it. This works fine for the continuation to $\Real(z)<1$, but for
$\Real(z)\geq 1$ there is a problem with applying $M$ (which we
must do before $Q_\lambda$) as the singularity at $0$ is of the
kind  $u^{-z}$ and is not integrable anymore. So we apply first
$L^2=(1-M)^2$ and only later $U$.

We compute:
\begin{eqnarray*}
M(t^{-z}\Un_{0<t<1})(u) &=& \frac{\int_0^{\min(1,u)}
t^{-z}\,dt}{u} = \frac{u^{-z}\Un_{u\leq1}(u)}{1-z}+
\frac{1}{1-z}\frac{\Un_{u>1}(u)}{u}\\ M^2(t^{-z}\Un_{0<t<1})(u)
&=& \frac{u^{-z}\Un_{u\leq1}(u)}{(1-z)^2}  +
\frac{1}{(1-z)^2}\frac{\Un_{u>1}(u) }{u}+
\frac{1}{1-z}\frac{\log(u)\Un_{u>1}(u)}{u}\\
L^2(t^{-z}\Un_{0<t<1})(u) &=&
\frac{z^2}{(z-1)^2}\,u^{-z}\Un_{u\leq1} + (\frac{z^2}{(z-1)^2} -
1)\frac{\Un_{u>1}(u)}{u} +
\frac{1}{1-z}\frac{\log(u)\Un_{u>1}(u)}{u}\\
\end{eqnarray*}

We note that $U(\Un_{u>1}/u) = UI(\Un_{u<1}) =
(M-1)\cF_+(\Un_{u<1}) = (M-1)(\sin(2\pi u)/(\pi u))$ is $O(1/u)$
(from the existence of the Dirichlet integral) for $u>\lambda$
and then that $U(\log(u)\Un_{u>1}/u) = UMI(\Un_{u<1}) =
MUI(\Un_{u<1}) = M(M-1) (\sin(2\pi u)/(\pi u)$ is
$O((1+|\log(u)|)/u)$. Clearly this reduces the problem of
$V(t^{-z}\Un_{0<t<1})(u)$ to the problem of the analytic
continuation and estimation of $U(t^{-z}\Un_{0<t<1})(u)$. From
\cite{Bal3}, proof of  Lemme 6, one has
$$U(t^{-z}\Un_{0<t<1})(u) = \frac{\sin(2\pi u)}{\pi u} +
\frac{z}{\pi u}\int_1^\infty t^{z-1}\sin(2\pi ut)\,\frac{dt}{t}$$
and (for example) from \cite{cras2} we know that the integral is
an entire function of $z$ which is $O(1/u)$ on
$[\lambda,\infty)$, uniformly in $z$ when $|z|$ is bounded. We
also see from this and from the integral representation of
$\widehat{B}(w)$ that it has at most a pole of order $2$ at $w=0$
(which is $z=1$).

The functional equation holds on the critical line from the
spectral representation of $\cG_\Lambda = C_\lambda I V$, hence
it holds on $\CC$ by analytic continuation. As $V(1-w)$ has poles
at $1-w = -2, -4, \dots$, and the left hand side is regular at
these values of $w$ it follows that  $\widehat{B}(1-w)$ has to
vanish for $1-w = -2, -4, \dots$.
\end{proof}

\begin{note}
The distance function $D(\lambda)^2$ has two components: one
corresponding to the distance to the subspace $M_\Lambda$ in
$\HH^2$ and then another one corresponding to the additional
distance inside this space to the translates $C_\theta(A)$,
$\lambda\leq\theta\leq1$. The first step has absolutely no
arithmetic, it is a problem of analysis. In the second step the
orthogonal projections to $M_\Lambda$ of the vectors
$Y^\lambda_{\rho,k}$, $\zeta(\rho) = 0$, $k<m_\rho$ are
obstructions.  When $\lambda\to0$ ($\Lambda\to\infty$) the first
contribution is presumably much smaller than the second, and the
original vectors $Y^\lambda_{\rho,k}$ will not themselves differ
much from their orthogonal projections to $M_\Lambda$. This seems
to suggest as a plausible thing that the estimate $(\sum_\rho
m_\rho^2/|\rho|^2)/|\log(\lambda)|$ gives the exact asymptotic
decrease of $D(\lambda)^2$ (under assumption of the Riemann
Hypothesis).
\end{note}

\section{Sonine spaces of de~Branges, novel spaces $H\!P_\lambda$,
  vectors $Z^\lambda_{\rho,k}$, Krein string of the zeta function}

Let $K$ be the Hilbert space $L^2((0,\infty),dt)$ of
complex-valued square-integrable functions on $(0,\infty)$ with
Hilbertian scalar product $(f,g) = \int_0^\infty
f(t)\overline{g(t)}\,dt$. We also use the ``Euclid'' bilinear
form $[f,g] = \int_0^\infty f(t)g(t)\,dt$. A vector $Z(t)$ is
``Euclid-perpendicular'' to a subspace $H$ for the bilinear form
$[f,g]$ if only and if $\overline{Z(t)}$ is
(``Hilbert'')-perpendicular to $H$ for the scalar product $(f,g)$
if and only if $Z(t)$ is Hilbert-perpendicular to the
complex-conjugated space $\overline H$. We also consider the
functions in $K$ as even functions with the definition $f(t) =
f(|t|)$ for $t<0$.

The Mellin transform (which is taken in the $L^2$-sense for
$\Real(s) = \frac12$)
$$f(t)\mapsto \wh f(s) = \int_0^\infty f(t) t^{-s}\,dt$$
isometrically identifies $K$ with the Hilbert space $L^2(s =
\frac12 + i\tau, d\tau/2\pi)$. The cosine transform $\cF_+$ acts
(in the $L^2$ sense) on $K$ as $\cF_+(f)(t)=2\int_0^\infty
\cos(2\pi tu)f(u)du$. It is a real operator. One has $\cF_+^2 =
1$, so $K$ is the orthogonal sum of the subspaces of invariant
functions (``self-reciprocal'') under $\cF_+$ and the subspaces
of anti-invariant (``skew-reciprocal'') functions.  The operator
$I$ is $f(t)\mapsto f(1/t)/|t|$. The composite $\Gamma_+ = \cF_+
I$ is scale invariant so it is diagonalized by the Mellin
transform: $\wh{\Gamma_+(f)}(s) = \chi_+(s) \wh f(s)$. This is
also written as
$$\wh{\cF_+(f)}(s) = \chi_+(s) \wh f(1-s)$$
The function $\chi_+(s)$ is a meromorphic function in the complex
plane which is related to the Tate Gamma function $\gamma_+(s)$
through $\chi_+(s)= \gamma_+(1-s) = \gamma_+(s)^{-1}$. One has:
$$\chi_+(s) = \pi^{s-1/2}\frac{\Gamma((1-s)/2)}{\Gamma(s/2)} =
2^{s}\pi^{s-1}\sin(\frac{\pi s}2)\Gamma(1-s) =
\frac{\zeta(s)}{\zeta(1-s)}$$
So an even function is self-reciprocal under the cosine transform
if and only if its \emph{right} Mellin transform satisfies the
zeta-functional equation. Under the \emph{left} Mellin transform
$f(t)\mapsto \int_0^\infty f(t) t^{s-1}\,dt$, which is the one
usually used in discussing the functional equation,
self-reciprocal functions under the cosine transform satisfy the
functional equation of $\zeta(1-s)$. The M\"untz formula
\cite[II.11]{Titchmarsh} shows that, for suitably regular
functions $f(t)$, the scale-invariant operator corresponding to
$\zeta(s)$ when using the \emph{left} Mellin transform is given
explicitely as a modified Poisson summation. So the action of the
scale invariant operator corresponding to $\zeta(s)$ under the
\emph{right} Mellin transform is expressed as the co-Poisson
summation (when applied to suitably regular functions; a more
detailed analysis will be given in \cite{copoisson}).

\let\la=\lambda \let\La=\Lambda

The intertwining property for co-Poisson summations ($t>0$):

$$\cF_+\left(\sum_{n\geq1} \frac{\alpha(t/n)}{n} - \int_0^\infty
\frac{\alpha(1/t)}{t}\,dt\right) = \sum_{n\geq1}
\frac{\alpha(n/t)}{t} - \int_0^\infty {\alpha(t)}\,dt$$

is an equivalent expression of the zeta-functional equation. It
shows how to give examples of even functions $f(t)$ which vanish
identically in a neighborhood $(-\la,\lambda)$ of the origin and
such that their Fourier cosine transform has the same
property. For this we take $\alpha(t)$, smooth with support in
$[\lambda,\Lambda]$ ($\Lambda = 1/\lambda$) and such that
$\int_0^\infty \alpha(t)\,dt = 0 = \int_0^\infty
(\alpha(1/t)/t)\,dt$. A non zero function $f(t)$  may be obtained
this way only for $0<\lambda<1$.

Nevertheless there exists for arbitrary $\lambda>0$ non-zero
square integrable even functions $f(t)$ vanishing in
$(0,\lambda)$, and such that $\cF_+(f)(t)$ also vanishes in
$(0,\lambda)$. To the best of the author's knowledge this was
first put forward as a fact of special importance in Analysis by
de~Branges in \cite{Bra64}. There, a beautiful isometric
expansion of self- and skew-reciprocal functions for the Hankel
transform of zeroth order is proven. The cosine transform is
(essentially) the Hankel transform of order $-\frac12$. The sine
transform is (essentially) the Hankel transform of order
$+\frac12$.

\begin{definition}
We let $K_\la\subset K$ be the Hilbert space of square-integrable
(even) functions $f(t)$ vanishing in $(0,\lambda)$, and such that
$\cF_+(f)(t)$ also vanishes in $(0,\lambda)$.
\end{definition}

An explicit example of a function having this property, but which
is not square-integrable, arises from an integral formula of
Sonine concerning Bessel functions \cite[p. 38]{Sonine}. This
example is in Titchmarsh's book on Fourier integrals
\cite[9.12.(8)]{TitchmarshFourier}. The analogous example which
is associated to the sine transform is a square-integrable (odd)
function. The proof given by de~Branges for existence of an even
square-integrable $f(t)$ with the Sonine property appears in
\cite{Bra64} on top of page 449. He constructs directly its 
Mellin transform from a trick which makes use of the already
known non-triviality of the spaces with the Sonine property for
the Hankel transform of positive order. Here is another trick: we
take a non-zero square-integrable odd function $g(t)$ which works
for the sine transform but with $\la^\prime = \la+1$. Let $f(t) =
g(t-1) - g(t+1)$. Then $f(t)$ is even, and non-trivial. It
vanishes on $(-\la,+\la)$ and its Fourier (cosine) transform
vanishes on $(-\la^\prime,+\la^\prime)$.

Actually the simplest method leading to explicit examples of
Sonine functions (in the Schwartz class), with arbitrarily large
$\lambda$, was communicated to the author by Professor
Kahane\footnote{\rm Letter to the author, March 22, 2002}:
\setcounter{footnote}0
the first observation is that it is enough to regularize a
tempered \emph{distribution} having the Sonine property to a
Schwartz function: additive convolution with a test-function
supported in $(-\epsilon,+\epsilon)$, and multiplication with its
Fourier (this replaces $\la$ with $\la-\epsilon$). The second is
that it is easy to obtain such distributions from the Poisson
distribution $\sum_{n\in\ZZ} \delta(x-n)$. As an example of how
to proceed one may take
$$x^3\prod_{1\leq j\leq N} (x^2 - j^2)^2\;\sum_{n\in\ZZ}
\delta^{\,\prime}(x-n)$$
and replace $x$ by $\sqrt N x$. This gives an even tempered
distribution with the Sonine property for $\la < \sqrt N$. In
\cite{copoisson} we use the \emph{multiplicative convolution} to
regularize such distributions and thus obtain more general
co-Poisson intertwining formulae.

An existence proof of Sonine square-integrable functions is
straightforward: it suffices to say that
$L^2(0,\lambda)+\cF_+(L^2(0,\lambda))$ is a closed (obviously
proper) subspace of $L^2(0,\infty)$. The (thus non trivial)
perpendicular complement is the  space $K_\la$. That the space
sum is closed follows readily  (see \cite[sect 2.9,
p. 126-127]{Dym1}) from the fact that the compact operator
$P_\la\cF_+ P_\la$ (where $P_\la$ is orthogonal projection to
$L^2(0,\lambda)$) has operator bound strictly less than one (as
no function can be compactly supported and with its Fourier
compactly supported). The next step is to actually write down
explicitely the associated orthogonal projection. This has led
the author recently to advances in the theory of the de~Branges
Sonine spaces (\cite{cras3}).

De~Branges proves that for each $f\in K_\la$ its completed
(right) Mellin Transform $$M(f)(s) := \pi^{-s/2}\Gamma(\frac
s2)\int_\la^\infty f(t)t^{-s}dt$$ is an entire function. The
proof (up to a change of variable) appears on page 447 of
\cite{Bra64}. We gave in \cite{cras2} another, more elementary,
proof.  The space of entire functions thus associated to $K_\la$
is among the ``Sonine spaces'' from \cite{Bra}, also studied in
\cite{Rov66, Rov67, Rov69}. The Sonine spaces are Hilbert spaces
of entire functions satisfying the axioms of \cite{Bra}. We will
also call $K_\la \subset L^2(0,\infty)$ a Sonine space.

\begin{note}
In \cite{Bra} as well as in other cited references it is the
horizontal axis which is the axis of symmetry. Comparison with
our conventions requires a change of variable (such as $s =
\frac12 - 2iz$), as it is the critical line which we use as the
axis of symmetry for the Hilbert spaces of entire functions.
\end{note}

So $K_\la = L^2((\la,\infty),dt)\cap
\cF_+(L^2((\la,\infty),dt))$. The convention in force in this
chapter will be to use the \emph{right} Mellin transform:
$$\widehat{f}(s) = \int_0^\infty f(t)t^{-s}\,dt$$

\begin{theorem}
The spaces $K_\la$ are all non-reduced to $\{0\}$. The Mellin
transforms of elements $f(t)$ from $K_\la$ are entire functions
with trivial zeros at $s=-2n$, $n\in\NN$. The entire functions
$$M(f)(s) = \pi^{-s/2}\Gamma(\frac s2)\widehat{f}(s)$$ satisfy
the functional equations
$$M(\cF_+(f))(s) = M(f)(1-s)$$
For each $w\in\CC$, each $k\in\NN$, the linear forms $f\mapsto
M(f)^{(k)}(w)$ are continuous and correspond to (unique) vectors
$Z^\la_{w,k}\in K_\la$: $\forall f\in K_\la\ [f, Z^\la_{w,k}] =
M(f)^{(k)}(w)$.
\end{theorem}

\begin{proof}
As we said, most of this is, up to a change of variable, from
\cite{Bra64}. The $L^2$-boundedness of the evaluations of the
derivatives ($k\geq1$) follow by the Banach-Steinhaus theorem
from the case $k=0$. We provide elementary proofs of all
statements in \cite{cras2}.
\end{proof}

\begin{note}
Evaluators such as $Z^\la_{w,k}$ for $k\geq1$, which are
associated to derivatives, do not seem to have been put to use so
far either in the general theory \cite{Bra}, or in the special
theory of Sonine spaces \cite{Bra64, Bra92, Bra94}.
\end{note}

\begin{note}
We have changed our conventions from  \cite{cras2} where we were
studying the functions in $L^2(0,\La)\cap I\cF_+I (L^2(0,\La)) =
I(K_\la)$ using left Mellin transforms. Here we study the
functions from $K_\la$ using the right Mellin transforms. So we
deal with exactly the same entire functions in the complex plane.
\end{note}

\begin{prop}
One has $K_\la = \cap_{\mu<\la} K_\mu$ and $K_\la =
\overline{\cup_{\mu>\la} K_\mu}$. Furthermore $L^2((0,\infty),dt)
= \overline{\cup_{\la>0} K_\la}$.
\end{prop}

\begin{proof}
Directly from the definition, the $K_\la$'s form a decreasing
chain as $\la\to\infty$ and the first statement holds. Let $\phi$
be perpendicular to each $K_\mu$ for $\mu>\la$. Then $\phi \in
L^2(0,\mu)+\cF_+L^2(0,\mu)$ and the decomposition as $f^\mu +
g^\mu$ is unique. The entire function $g^\mu$ must then not
depend on $\mu>\la$, as the difference between two such will be
compactly supported, hence zero. Then $f^\mu$ does not depend on
$\mu$ either and is in $L^2(0,\la)$ and $g^\mu\in
\cF_+L^2(0,\la)$. So $K_\la = \overline{\cup_{\mu>\la}
K_\mu}$. The same proof shows the last statement.
\end{proof}

\begin{prop}
One has $\cF_+(Z^\la_{w,k}) = (-1)^k Z^\la_{1-w,k}$.
\end{prop}

\begin{proof}
Using the Euclid bilinear form: $[\cF_+(Z^\la_{w,k}),f] =
[Z^\la_{w,k},\cF_+(f)] = M(\cF_+(f))^{(k)}(w) = (-1)^k
M(f)^{(k)}(1-w)$ from $M(\cF_+(f))(w) = M(f)(1-w)$.
\end{proof}

We use the bilinear form $[f, Z^\la_{w,k}]$ so that the vectors
$Z^\la_{w,k}$ depend analytically on $w$. Evaluators ($k=0$)
associed to points off the symmetry axis are always non-zero
vectors in de~Branges spaces. Using our elementary techniques we
proved a stronger statement in the case at hand:

\begin{theorem}[\cite{cras2}]
Any finite collection of vectors $Z^\la_{w,k}$ is a linearly
independent system. In particular the vectors $Z^\la_{w,k}$ are
all non-vanishing.
\end{theorem}

\begin{note}
If we take an arbitrary sequence of distinct complex numbers
having an accumulation point the corresponding evaluators
$Z^\la_{w,0}$ will span $K_\la$. Remarkable orthogonal bases
consisting of evaluators $Z^\la_{w,0}$ exist as a general fact from
\cite{Bra}. Some other non-trivial examples of infinite and
minimal collection of evaluators are also known \cite{twosys}.
\end{note}

It is useful to have at our disposal ``augmented Sonine'' spaces
$L_\la\supset K_\la$,  whose elements' (Gamma-completed) Mellin
Transforms may have poles at $0$ and $1$. Let $N$ be the unitary
invariant operator which, under the right Mellin transform, has
spectral function $s/(s-1)$. Explicitely:

$$N(f)(t) = f(t) - \int_t^\infty \frac{f(u)}u\,du$$

Let $L_\la = N\cdot L^2((\la,\infty),dt)\cap \cF_+\cdot N\cdot
L^2((\la,\infty),dt)$.

\begin{theorem}
Let $\la>0$. Let $f\in L_\la$. The Mellin transform $\wh{f}(s) =
\int_0^\infty f(t)t^{-s}\,dt$ is an analytic function in
$\CC\setminus\{1\}$ with at most a pole of order $1$ at $s=1$. It
has trivial zeros at $s=-2n$, $n\geq1$. The function $M(f)(s) =
\pi^{-s/2}\Gamma(s/2)\widehat{f}(s)$, analytic in
$\CC\setminus\{0, 1\}$,  satisfies the functional equation
$$M(\cF_+(f))(s) = M(f)(1-s)$$
For each $w\in\CC\setminus\{0, 1\}$, each $k\in\NN$, the linear
forms $f\mapsto M(f)^{(k)}(w)$ are continuous.
\end{theorem}

\begin{proof}
As an intersection $L_\la$ is a closed subspace of
$L^2((0,\infty),dt)$ hence a Hilbert space. The square integrable
function $f(t)$ is a constant $\alpha(f)$ for $0<t<\la$ (which is
a continuous linear form in $f$). So $\int_0^\infty
f(t)t^{-s}\,dt$ is absolutely convergent and analytic at least
for $1/2<\Real(s)<1$. In this strip we may write it as:
$$\int_0^\infty f(t)t^{-s}\,dt = \frac{\alpha(f)\la^{1-s}}{1-s} +
\int_\la^\infty f(t)t^{-s}\,dt$$
which gives its analytic continuation to the right half-plane
$\Real(s)>1/2$ with  at most a pole at $s=1$. Let us also note
that the evaluation at these points are clearly continuous for
the Hilbert structure. We have:
$$\int_\la^\infty f(t)t^{-s}\,dt = \int_0^\infty
\cF_+(f)(u)\cF_+(\Un_{t>\la}t^{-s})(u)\,du$$
We known from \cite[Lemme 1.3.]{cras2} that the function
$\cF_+(\Un_{t>\la}t^{-s})(u)$ is an entire function of $s$, which
is (uniformly for $|s|$  bounded) $O(1/u)$ on $(\la,\infty)$, and
also that it is $\chi_+(s)u^{s-1} + O(1)$ on $(0,\la)$ (uniformly
for $\Real(s)\leq1-\epsilon<1$). Moreover $\cF_+(f)(u)$ is a
constant in the interval $(0,\la)$ (from $f\in L_\la$).
Combining these informations we get that 
the above displayed equation has an analytic continuation to the
critical strip $0<\Real(s)<1$. In this critical strip we have the
functional equation:
$$\wh{f}(s) = \chi_+(s)\wh{\cF_+(f)}(1-s)$$
as it holds on the critical line.  From this we get the
analytic continuation of $\wh f(s)$ to $\Real(s)<1$. We note that
$\chi_+(s)$ vanishes at $s=0$ and that this counterbalances the
(possible) pole of $\wh{\cF_+(f)}(1-s)$. Also this functional
equation shows that $\wh{f}(s)$ vanishes at $s = -2n$,
$n\geq1$. The evaluations at points strictly to the right of the
critical line are continuous, hence also at points to the left,
hence everywhere (except  of course at $s=0, s=1$, where instead
one may consider the residues) from the Banach-Steinhaus theorem.
\end{proof}

There are (for $w\neq 0,1$) evaluators $W^{\la}_{w,k} \in L_\la$
which project orthogonally to the evaluators $Z^\la_{w,k} \in
K_\la$. The augmented Sonine spaces $L_\la \supset K_\la$ are
natural for discussing properties of the zeta-function along the
lines involving the co-Poisson formula. Here we will stay in the
realm of the spaces $K_\la$, using the spaces $L_\la$ as an
auxiliary help.

Let $0<\la<1$ and let $\La = 1/\la$. Let $\alpha(t)$ be a smooth
function with support in $[\la,\La]$ and let $f(t)$ be the
co-Poisson summation $\sum_{n\geq1} \alpha(t/n)/n - \int_\la^\La
\alpha(t)\,dt/t$. The function $f(t)$ is a Schwartz function,
hence square-integrable. From the co-Poisson formula it belongs
to $L_\la$. If we impose the conditions that $\wh \alpha(0) = 0 =
\wh \alpha(1)$ then $f(t)$ belongs to the Sonine space
$K_\la$. At the level of Mellin transform, we have $\wh f(s) =
\zeta(s) \wh \alpha(s)$. So $f(t)$ is (Euclid)-perpendicular to
the evaluators $Z^{\la}_{\rho,k}$, $k<m_\rho$ associated with the
non-trivial zeros of the Riemann zeta function and with their
(eventual) multiplicities. And conversely it follows
(\cite{cras2}) from $\wh f(s) = \zeta(s) \wh \alpha(s)$  that an
evaluator $Z^{\la}_{w,k}$ is (Euclid or Hilbert) perpendicular
to all functions $\sum_{n\geq1} \alpha(t/n)/n$ with $\alpha$
smooth function with support in $[\la,\La]$ and  $\wh \alpha(0) =
0 = \wh \alpha(1)$ if and only if $w$ is a non-trivial zero of
the zeta function with multiplicity strictly bigger than $k$.

\begin{definition}
Let $0<\la<1$ and $\La = 1/\la$. We let $W_\la$ be the closure in
$K$ of the functions $\sum_{n\geq1} \alpha(t/n)/n$, with
$\alpha(t)$ \emph{smooth} with support in $[\la,\La]$, and such
that $\widehat{\alpha}(0) = 0 = \widehat{\alpha}(1) $.
\end{definition}

\begin{definition}
Let $0<\la<1$ and $\La = 1/\la$. We let $W_\la^\prime$ be the
sub-vector space of $K$ comprising the square-integrable
functions $f(t)$ which may be written as $\sum_{n\geq1}
\alpha(t/n)/n$, where $\alpha(t)\in L^1(\la,\La)$ and
$\widehat{\alpha}(0) = 0 = \widehat{\alpha}(1) $.
\end{definition}

\begin{definition}
Let $0<\la<\infty$. We let $Z_\la$ be the closed subspace of
$K_\la\subset K$ spanned by the evaluators $Z^{\la}_{\rho,k}$,
$0\leq k<m_\rho$ associated with the non-trivial zeros of the
Riemann zeta function and with their (eventual) multiplicities.
\end{definition}

The main theorem (whose proof takes up the next pages) is:

\begin{theorem}\label{maintheorem}
{\bf 1.} Let $0<\la<1$. One has $W_\la\subset W_\la^\prime
\subset K_\la$. The subspace $W_\la^\prime$ is closed and equals
$\cap_{0<\mu<\la} W_\mu = \cap_{0<\mu<\la} W_\mu^\prime$. One has
$W_\la = \overline{\cup_{\la<\mu<1} W_\mu} =
\overline{\cup_{\la<\mu<1} W_\mu^\prime}$. One has $K_\la =
W_\la^\prime\perp Z_\la$.\\ {\bf 2.} The set of $\la$'s for which
$W_\la\subset W_\la^\prime$ is a strict inclusion is at most
countable.\\ {\bf 3.} Let $1\leq \la <\infty$. One has $K_\la =
Z_\la$.
\end{theorem}

\begin{definition}
Let $0<\la<1$. We let $H\!P_\la$ be the perpendicular complement
in $K_\la$ of $W_\la$.
\end{definition}

We thus have $H\!P_\la \supset Z_\la$ and the question whether
this may be strict is interesting (equivalently whether $W_\la
\subset W_\la^\prime$ may be a strict inclusion). This question
is related to the properties of the Krein spaces of entire
functions of finite exponential type which are  associated with
the measure $|\zeta(\frac12 + i\tau)|^2 d\tau/2\pi$ on the
critical line. From $W_\la^\prime = \cap_{\mu<\la} W_\mu$ a
strict inclusion may happen only for a countable set of $\la$'s.

As usual our axis of symmetry is the critical line, not the real
axis, and we use the Mellin transform to define Paley-Wiener
functions, not the additive Fourier transform. Let $\mu$ be a
measure on the critical line, and let $H = L^2(s = \frac12 +
i\tau, d\mu)$. We suppose $1/s \in H$ and $d\mu(\frac12 + i\tau)
= d\mu(\frac12 - i\tau)$. Let $\La\geq1$ and let $I^\La$ be the
subspace of $H$ of ($\mu$-equivalence classes of) functions
$F(s)$ which are also entire functions of exponential type at
most $\log(\La)$. Let $J^\La$ be the subspace of $H$ of functions
$F(s)$ which are also entire functions of exponential type
strictly less than $\log(\La)$ (for $\La=1$ this means $J^1 =
\{0\}$). 
It is proven in \cite{Dym} that $I^\La$, if it does not span
$H$, is a closed subspace. It will then contain the closure
of $J^\La$ and the question whether it may be strictly
larger is subtle. 
An isometric representation exists, the \emph{Krein
string}, where, if the description of the string is complete
enough, one may read the answer to the question. We do not go
into more details and refer the reader to the book \cite{Dym}
which is devoted to the theory of the Krein string, and which
also contains an introduction to the de~Branges theory. The
following theorem is due to Krein and is also fundamental in the
general de~Branges theory.

\begin{theorem}[Krein, \cite{Kre}]\label{kreintheorem}
Let $F(z)$ be an entire function which is in the Nevanlinna class
separately in the half-plane $\Imag(z)>0$ and in the half-plane
$\Imag(z)<0$. Then $F(z)$ has finite exponential type which is
given by the formula
$$\max(\limsup_{\sigma\to+\infty}\frac{\log|F(i\sigma)|}{\sigma},
\limsup_{\sigma\to+\infty}\frac{\log|F(-i\sigma)|}{\sigma})$$
\end{theorem}

We recall that one possible definition of the Nevanlinna class of
a half-plane is as the space of quotients of bounded analytic
functions. For example it is known that any function in the Hardy
space of a half-plane is a Nevanlinna function. Krein's theorem
is more complete but we only need the result given here. Of
course we will be using this theorem with the critical line
replacing the horizontal axis.

For the following steps  we let $0<\la\leq1$, $\La = 1/\la$, and
the notations $H$, $I^\La$, $J^\La$ are relative to the measure
$d\mu(s) = |\zeta(s)|^2 d\tau/2\pi$ on the critical line ($s =
\frac12 + i \tau$).
$$H = L^2(\Real(s) = \frac12, |\zeta(s)|^2 \frac{d\tau}{2\pi})$$
Unfortunately we are unable to describe the associated Krein
string. Rather we will explain how to isometrically identify the
co-Poisson spaces $W_\la$ (resp. $W_\la^\prime$) (here $0<\la<1$)
with subspaces of codimension $2$ of $\overline{J^\La}$
(resp. $I^\La$). This will be used in the proof of the main
theorem \ref{maintheorem}.

\begin{lem}
A function $G\in H$ is perpendicular to $J^\La$ if and only if it
is perpendicular to all functions $(u^s -1)/s$ for $\la\leq u\leq
\La$. Hence the closure $\overline{J^\La}$ is also the closure of
the finite linear combinations $(u^s -1)/s$ for $\la\leq u\leq
\La$.
\end{lem}

\begin{proof}
One direction is obvious. Let us now assume that $G\perp (u^s
-1)/s$ for $\la\leq u\leq \La$. Let $F\in J^\La$ and let
$\epsilon>0$ be such that the type of $F$ is
$<\log(\La)-\epsilon$. We consider
$$\int F(s)\frac{e^{\epsilon s} -
1}{s}\overline{G(s)}|\zeta(s)|^2\,d\tau$$
If we take $F(s) = u^s$ with $e^\epsilon\la\leq u\leq
e^{-\epsilon}\La$ this integral vanishes. Using the
Pollard-de~Branges-Pitt ``lemma'' (\textit{sic}) from \cite[4.8.,
p.108]{Dym}, we deduce that the integral with the original $F(s)$
vanishes too. Then from $|(e^{\epsilon s} - 1)/\epsilon s|\leq 2
(e^{\epsilon/2} - 1)/\epsilon$ and dominated convergence we get
the desired conclusion.
\end{proof}

\begin{lem}
Let $F(s)\in I^\La$. One has $F(s)\zeta(s)\in N\cdot\La^s\HH^2$
and also $F(1-s)\zeta(s)\in N\cdot\La^s\HH^2$ (we write $\HH^2$
for the Hardy space of the right half-plane and we recall that
$N$ is the operator of multiplication with $s/(s-1)$.)
\end{lem}

\begin{proof}
The product  $F(s)\zeta(s)$ belongs to $L^2(\Real(s)=1/2,
d\tau/2\pi)$.  As $I^\La \subset J^{\La \exp(\epsilon)}$ for
$\epsilon>0$, $F(s)\zeta(s)$ is in the closure of finite sums of
functions $(u^s -1)\zeta(s)/s$ for $e^{-\epsilon}\la\leq u\leq
e^{+\epsilon}\La$. It belongs to the closed space $N\cdot
(e^\epsilon\La)^s\HH^2$ as $\zeta(s)/s$ itself belongs to
$N\cdot\HH^2$. We note that this space is the image under $N$ of
the Mellin transform of $L^2((e^{-\epsilon}\la, \infty),dt)$ so
after letting $\epsilon\to0$ we obtain that $F(s)\zeta(s)$
belongs to $N\cdot \La^s\HH^2$.  We note that $F(s)\to F(1-s)$ is
an isometry of $I^\La$ and the conclusion then follows.
\end{proof}

\begin{theorem}
An entire function $F(s)$ belongs to $I^\La$ (i.e. it is in $H$
and of exponential type at most $\log(\La)$) if and only if
$F(s)\zeta(s)$ is the Mellin transform of an element in
$L_\la$. The space $I^\La$ is a closed subspace of $H$ and is
isometric through $F(s)\to\zeta(s)F(s)$ to the subspace of
$L_\la$ of functions whose Mellin transform vanish at the zeros
of the zeta function with at least the same multiplicities. For
each complex number $w$ the evaluations $F\mapsto F(w)$ are
continuous linear forms on $I^\La$.
\end{theorem}

\begin{proof}
\From\ the lemma $\zeta(s)F(s)$ is the Mellin transform of an
element of $N\cdot L^2((\la,\infty),dt)$ whose image under
$\cF_+$ also belongs to $N\cdot L^2((\la,\infty),dt)$ (as this
corresponds to the replacement $F(s)\mapsto F(1-s)$). So the map
$F(s)\to\zeta(s)F(s)$ is an isometric embedding into
$\wh{L_\la}$. If an element $G(s)$ from $\wh{L_\la}$ vanishes at
the non-trivial zeros of the zeta function (taking into account
the multiplicities) then it factorizes as $G(s) = F(s)\zeta(s)$
with an entire function $F(s)$ (as $G(s)$ also vanishes at the
trivial zeros and has at most a pole of order $1$ at
$s=1$).  From this, $F(s)$ is in the right half-plane in the
Nevanlinna class (of quotients of bounded analytic functions)
because both $F(s)\zeta(s)$ and $\zeta(s)$ are meromorphic
functions in this class. And the same holds in the left
half-plane, as $\cF_+(G)(s) = F(1-s)\zeta(s)$. We now use the
theorem of Krein \ref{kreintheorem} which tells us that the
entire function $F(s)$ has finite exponential type given by
$$\max(\limsup_{\sigma\to+\infty}\frac{\log|F(\sigma)|}{\sigma},
\limsup_{\sigma\to+\infty}\frac{\log|F(1-\sigma)|}{\sigma})$$
From  this formula, and from $F(s)\zeta(s)\in
N\cdot\La^s\HH^2$, $F(1-s)\zeta(s)\in N\cdot\La^s\HH^2$, and from
the fact that elements of $\HH^2$ are bounded in $\Real(s)\geq1/2
+ \epsilon > 1/2$, we deduce that the exponential type of $F(s)$
is at most $\log(\La)$. So $I^\La$ is isometrically identified
with the functions in $\wh{L_\la}$ vanishing at least as
$\zeta(s)$ does. This space is closed because the evaluators are
continuous linear forms on $L_\la$. From this we see that the
evaluators $F\mapsto F(s)$  are continuous linear forms except
possibly at the zeros and poles of $\zeta(s)$, and the final
statement then  follows from this and the Banach-Steinhaus
theorem (as $I^\La$ is a Hilbert space from the preceding).
\end{proof}

\begin{theorem}\label{omegaprime1}
Let $0<\la<\infty$.\\ 1. The vectors $Z^\la_{\rho,k}$,
$k<m_\rho$, span $K_\la$ if and only if $\la\geq1$.\\ 2. A
function $\alpha(s)$ on $\Real(s)=1/2$ is the Mellin transform of
an element of $K_\la$ perpendicular to $Z_\la$ if and only if:\\
\indent\hskip1cm a. It is square integrable on the critical line
for $d\tau/2\pi$.\\ \indent\hskip1cm b. One has
$\alpha(s)=\zeta(s)s(s-1)\beta(s)$ with $\beta(s)$ an entire
function of finite exponential type at most $\log(1/\la)$.
\end{theorem}

\begin{proof}
\From\ the existence of $W_\la$ the vectors $Z^\la_{\rho,k}$,
$k<m_\rho$, do not span $K_\la$ if $\la<1$. Let $f\in
(Z_\la)^\perp \cap K_\la$. By definition its Mellin transform
vanishes at the non-trivial zeros of $\zeta$. It also vanishes at
the trivial zeros and at $0$ so it may be written
$$\wh{f}(s) = s(s-1)\zeta(s)\theta(s)$$
with an entire function $\theta(s)$. In the right half-plane
$\theta(s)$ is in the Nevanlinna class (of quotients of bounded
analytic functions) because both $\wh{f}(s)$ and $s(s-1)\zeta(s)$
are meromorphic in this class. From  the functional equation one
has
$$\wh{\cF_+(f)}(s) = s(s-1)\zeta(s)\theta(1-s)$$
So $\theta(s)$ is in the Nevanlinna class of the left
half-plane. We now use the theorem of Krein \ref{kreintheorem}
and conclude that the entire function $\theta(s)$ has finite
exponential type which is given as
$$\max(\limsup_{\sigma\to+\infty}\frac{\log|\theta(\sigma)|}{\sigma},
\limsup_{\sigma\to+\infty}\frac{\log|\theta(1-\sigma)|}{\sigma})$$
This formula (elements of $\HH^2$ are bounded in $\Real(s)\geq1$)
shows that the exponential type of $\theta(s)$ is at most
$\log(1/\la)$. This shows $Z_\la = K_\la$ for $\la>1$. Let us
prove this also for $\la=1$: on the line $\Real(s) = +2$ one has
$\wh{f}(s) = O(1)$ (as it belongs to $\HH^2(\Real(s)>1/2)$) hence
$\theta(s)$ is $O(1/s(s-1))$. So it is square integrable on this
line and by the Paley-Wiener theorem it vanishes identically as
it is of minimal exponential type.

Conversely, let $F(s) = s(s-1)\beta(s)$ be an entire function of
finite exponential type at most $\log(1/\la)$ which is such that
$\alpha(s) = \zeta(s)F(s)$ is square-integrable on the critical
line. From the previous theorem $F(s)$ is in the closed subspace
$I^\La$ of $H$ and $\alpha(s)$ is the Mellin transform of an
element $f(t)$ of $L_\la$. As $\alpha(s)$ is analytic at $s=1$
and vanishes at $s=0$ one has in fact $f\in K_\la$. And
$\alpha(s)$ vanishes at the zeros of zeta with at least the same
multiplicities, in other words $f$ is perpendicular to $Z_\la$.
\end{proof}

\begin{lem}\label{lemmeO1}
Any function $F$ in $I^\La$ is $O(1)$ in the closed strip $-1\leq
\Real(s)\leq2$, in particular on the critical line.
\end{lem}

\begin{proof}
\From\  the fact that $F(s)\zeta(s)\frac{s-1}{s}\La^{-s}$ is
bounded on the line $\Real(s) = 2$ (as it belongs to the Hardy
space $\HH^2(\Real(s)>\frac12)$) one deduces that $F(s)$ is
bounded on $\Real(s) = 2$ , hence also on $\Real(s) = -1$ (as
$F(1-s)$ also belongs to $I^\La$.) As it has finite exponential
type we may apply the Phragmen-Lindel\"of theorem to deduce  that
$F(s)$ is $O(1)$ on this closed vertical strip.
\end{proof}

\begin{lem}
Let $\La >1$. Let $K^\La$ be the closure of $J^\La$ in $H$. Let
$K^\La_0$ be the subspace of functions in $K^\La$ vanishing at
$0$ and at $1$, and similarly let $J^\La_0$ be the subspace of
$J^\La$ of functions vanishing at $0$ and at $1$. Then $K^\La_0$
is the closure of $J^\La_0$.
\end{lem}

\begin{proof}
Let for $1<\mu<\La$:
$$A_\mu(s) = \frac{(\mu^{s/2} - 1)(\mu^{s/2} -
\mu^{1/2})}{\log(\mu)(1-\mu^{1/2})/2}\frac{1}{s}$$
This is an entire function of exponential type $\log(\mu)$, in
$H$ and with $A_\mu(0) = 1$, $A_\mu(1) = 0$. Let also $B_\mu(s) =
A_\mu(1-s)$. Let $F\in K^\La_0$ and let us write $F=\lim F_\mu$
with $F_\mu\in J^\mu$, $\mu<\La$. One has $F_\mu(0)\to F(0) = 0$
and $F_\mu(1)\to F(1) = 0$, because evaluations are continuous
linear forms on $I^\La$. So $F = \lim (F_\mu - F_\mu(0)A_\mu -
F_\mu(1)B_\mu)$ (clearly the norms of $A_\mu$ and $B_\mu$ are
bounded as $\mu\to\La$).
\end{proof}

\begin{theorem}\label{omega}
Let $0<\la<1$. A function $\alpha(s)$ on $\Real(s)=1/2$ is the
Mellin transform of an element of $W_\la$ if and only if:\\ 1. It
is square integrable on the critical line for $d\tau/2\pi$.\\
2. It is in the closure of the square integrable functions
$\alpha(s)=\zeta(s)s(s-1)\beta(s)$ with $\beta(s)$ an entire
function of finite exponential type \emph{strictly less} than
$\log(1/\la)$.
\end{theorem}

\begin{proof}
Let $0<\la<1$. We have to show $\wh{W_\la} =
K^\La_0\cdot\zeta(s)$. First let us prove the inclusion
$\wh{W_\la}\subset K^\La_0\cdot\zeta(s)$: let $\phi(u)$ be a
smooth function with support in $[\la,\La]$ with $\wh{\phi}(1) =
0 = \wh{\phi}(0)$. It is elementary that there exists $\psi(u)$
smooth with support in $[\la,\La]$ and with $\wh{\phi}(s) =
s(s-1)\wh{\psi}(s)$. If we now consider for $a\to1^-$ the smooth
functions $\phi_a(u)$ with support in $[\la^a,\La^a]$ such that
$\wh{\phi_a}(s) = s(s-1)\wh{\psi}(as)$ then
$\wh{\phi_a}(s)\zeta(s)$ belongs to $K^\La_0\cdot\zeta(s)$ and
converge  to $\wh{\phi}(s)\zeta(s)$ in $L^2$ norm on the critical
line as $a\to1^-$.

For the converse inclusion $K^\La_0\cdot\zeta(s)\subset
\wh{W_\la}$ let $F\in K^\La_0$.  We may approximate $F$ with an
element of $J^\La_0$, so we may assume $F$ itself to be of
positive exponential type $\log(\mu)<\log(\La)$. Let $\theta$ be
a smooth function with support in $[1/e, e]$, with
$\wh{\theta}(1/2) = 1$. Let $\wh{\theta_\epsilon}(s) =
\wh\theta(\epsilon(s-1/2) + 1/2)$. Let $F_\epsilon =
\wh\theta_\epsilon F$. In $H$ the functions $F_\epsilon$ converge
to $F$. From \ref{lemmeO1} we know that $F$ is $O(1)$ on the
critical line so the functions $F_\epsilon$ are $O(|s|^{-N})$ for
any $N\in\NN$. From the Paley-Wiener theorem they are the Mellin
transforms of $L^2$ functions $f_\epsilon(t)$ with support in
$[e^{-\epsilon}\mu^{-1}, e^{\epsilon}\mu]$. For $\epsilon$ small
enough this will be included in $[\la,\La]$. From  the decrease
on the critical line the functions $f_\epsilon(t)$ are smooth. As
$\wh{f_\epsilon}(0) = 0 = \wh{f_\epsilon}(1)$ this tells us that
$F_\epsilon(s)\zeta(s)$ is the Mellin transform of a co-Poisson
summation of a \emph{smooth} function, and this implies that
$F(s)\zeta(s)$ belongs to the (Mellin transform of) $W_\la$, as
$W_\la$ is defined as the closure of the co-Poisson summations of
smooth functions whose Mellin transforms vanish at $0$ and at $1$.
\end{proof}

\begin{theorem}\label{thml2}
Let $F(s)$ be an entire function of finite exponential type. Then
$$\int_{\Real(s) = 1/2} |F(s)|^2|\zeta(s)|^2|ds|<\infty \
\Longrightarrow\ \int_{\Real(s) = 1/2} |F(s)|^2 |ds| <\infty$$
\end{theorem}

\begin{proof}
We want to prove that any function $F(s)$ in $I^\La$ is
square-integrable for the Lebesgue measure on the critical
line. We know from \ref{lemmeO1} that it is $O(1)$ in the closed
strip $-1\leq\Real(s)\leq2$. From this, if for $s$ in this open
strip we express $F(s)$ as a Cauchy integral with contributions
from the two vertical sides and two horizontal segments, the
contribution of the horizontal segments will vanish when they go
to infinity. So:
$$F(s) = \int_{{\bf Re}(s)=2}{F(z)\over z-s}{|dz|\over2\pi} -
\int_{{\bf Re}(s)=-1}{F(z)\over z-s}{|dz|\over2\pi}$$ On
$\Real(s)=2$, $F(s)$ is square integrable because
$F(s)\zeta(s)\frac{s-1}{s}\La^{-s}$ is, as it belongs to
$\HH^2(\Real(s)>\frac12)$. It is an important fact that Cauchy
integrals of $L^2$ functions on vertical line realize the
orthogonal projection to the Hardy space of the corresponding
half-plane. Hence the first integral above defines a function
square-integrable on each vertical line $\Real(s)<2$. And the
second integral similarly for $\Real(s)>-1$ ($F(1-s)$ satisfies
the same hypotheses as $F(s)$). So $F(s)$ is square-integrable on
the critical line (and in fact on each vertical line in the
complex plane.)
\end{proof}

The next theorem establishes   $K_\la = W_\la^\prime\perp Z_\la$:

\begin{theorem}\label{omegaprime2}
Let $\la<1$. The functions $A(u)$ in $Z_\la^\perp\ \cap K_\la$
are exactly the square-integrable functions which may be written
$\sum_{n\geq1} g(u/n)/n$, with an integrable function $g(u)$
supported in $[\la, \La]$ and such that $\int_0^\infty g(u)du =
\wh{g}(0) = 0$. The function $g(u)$ is necessarily
square-integrable and necessarily satisfies $\wh{g}(1) =
\int_0^\infty \frac{g(u)}u du = 0$.
\end{theorem}

\begin{note}
By a variant on the M\oe bius inversion formula from $A(u) =
\sum_{n\geq1} g(u/n)/n$ one has $g(u) = \sum_{n\geq1}
\mu(m)A(u/m)/m$ (and this is a finite sum for each $u>0$) in case
$A$ (hence $g$ and conversely) has support in $(\la,\infty)$. It
involves then in a neighborhood of each $u>0$ only finitely many
terms. If $g(u)$ has support in $[\la, \La]$ we can express it on
this interval as a finite combination of $A(u/m)/m$'s. So if $A$
is $L^2$ then $g$ had to be $L^2$ to start with. Also  we will
see that if $A$ is $L^2$ then $\int_0^\infty \frac{g(u)}udu$
necessarily vanishes.
\end{note}

\begin{proof}
Let $A\in Z_\la^\perp\ \cap K_\la$ and $\alpha(s) =
\wh{A}(s)$. We know that $\alpha(s) = \zeta(s)F(s)$ with $F(s)$
an entire function vanishing at $0$ and $1$ and of exponential
type at most $\log(\La)$. From \ref{thml2} we know that $F(s)$
is square-integrable on the critical line for the Lebesgue
measure. So the Paley-Wiener theorem implies $F(s) = \wh{g}(s)$
with $g(u)\in L^2([\la, \La])$. We have our function $g(u)$ in
$L^2(\la,\La)$ and we want to show that $A(u)$ is equal to $B(u)
= \sum_{n\geq1} g(u/n)/n$.  From  Fubini $\int_\la^\infty
B(u)u^{-s}du = \zeta(s)\wh{g}(s) = \wh{A}(s) = \int_\la^\infty
A(u)u^{-s}du$ for $\Real(s)>1$ and so $B(u) = A(u)$ (almost
everywhere from the unicity theorem for Fourier transforms of
$L^1$-functions). We have shown that each $A\in Z_\la^\perp\ \cap
K_\la$ may be written (uniquely) as $\sum_{n\geq1} g(u/n)/n$ with
$g\in L^2((\la, \La),du)$, $\wh{g}(0) = \wh{g}(1) = 0$. So it
belongs to $W_\la^\prime$.

For the converse let $g\in L^2((\la, \La),du)$ be such that $A(u)
= \sum_{n\geq1} g(u/n)/n$ is square-integrable. Its distribution
theoretic Fourier transform is (from \ref{copoissongeneral}):
$$\cF\left(\sum_{n\geq1} g(u/n)/n \right) = \sum_{n\geq1}
g(n/u)/u - \wh{g}(0) + \wh{g}(1)\delta_0$$ This distribution must
coincide with the function which is the $L^2$-Fourier transform
of $A(u)$ and so the square-integrability of $A(u)$ implies the
vanishing of $\wh{g}(1)$.

If we impose  $ \wh{g}(0) = 0$ the  Fourier transform of $A(u)$
is the function $\sum_{n\geq1} g(n/u)/u$ which again vanishes on
$(0,\la)$. So $A$ belongs to $K_\la$. Its Mellin transform is an
entire function which by Fubini for $\Real(s)>1$ equals
$\zeta(s)\wh{g}(s)$ hence also everywhere. The vector $A(u)$ is
thus perpendicular to the vectors $Z^\la_{\rho,k}$, which means
that $A\in Z_\la^\perp$. This completes the proof of $K_\la =
W_\la^\prime \perp Z_\la$.
\end{proof}

We also take note of:

\begin{prop}
The map ``${\cdot\over\zeta(s)}$'' from $W_\la^\prime$ to
$L^2(\la,\La)$ is bounded.
\end{prop}

\begin{proof}
Each $g(u)$ is expressed (on $(\la,\La)$) as a finite M\oe bius
sum in terms of the $A(u/m)/m$'s, with a number of summands
independent of $A$.
\end{proof}

The main theorem sums up almost everything that preceded:

\begin{theorem}[{\ref{maintheorem}}]
{\bf 1.} Let $0<\la<1$. One has $W_\la\subset W_\la^\prime
\subset K_\la$. The subspace $W_\la^\prime$ is closed and equals
$\cap_{0<\mu<\la} W_\mu = \cap_{0<\mu<\la} W_\mu^\prime$. One has
$W_\la = \overline{\cup_{\la<\mu<1} W_\mu} =
\overline{\cup_{\la<\mu<1} W_\mu^\prime}$. One has $K_\la =
W_\la^\prime\perp Z_\la$.\\ {\bf 2.} The set of $\la$'s for which
$W_\la\subset W_\la^\prime$ is a strict inclusion is at most
countable.\\ {\bf 3.} Let $1\leq \la <\infty$. One has $K_\la =
Z_\la$.
\end{theorem}

\begin{proof}
The basic inclusions $W_\la\subset W_\la^\prime \subset K_\la$
are a corollary to the co-Poisson intertwining formula. One has
$K_\la = Z_\la$ for $\la\geq1$ from Theorem
\ref{omegaprime1}. Let $0<\la<1$. From Theorem \ref{omegaprime2}
we have identified $W_\la^\prime$ as the perpendicular component
in $K_\la$ of $Z_\la$. From Theorem \ref{omegaprime1}
$W_\la^\prime$ is isometrically identified with the closed
subspace of $L^2(\Real(s)=\frac12, |\zeta(s)|^2 d\tau/2\pi)$ of
(restrictions) of entire functions $F(s)$ of exponential type at
most $\log(1/\la)$ and vanishing at $0$ and at $1$. Hence
$W_\la^\prime = \cap_{0<\mu<\la} W_\mu^\prime$. From Theorem
\ref{omega} $W_\la$ is isometrically identified with the closure
in $L^2(\Real(s)=\frac12, |\zeta(s)|^2 d\tau/2\pi)$ of entire
functions $F(s)$ of exponential type \emph{strictly less than}
$\log(1/\la)$ and vanishing at $0$ and at $1$. Hence $W_\la =
\overline{\cup_{\la<\mu<1} W_\mu}$. Also $\la<\mu<1 \Rightarrow
W_\la\supset W_\mu^\prime \supset W_\mu$ (the last inclusion as
$W_\mu^\prime$ is known to be closed). Hence $W_\la =
\overline{\cup_{\la<\mu<1} W_\mu^\prime}$. Also $W_\la^\prime
\subset \cap_{0<\mu<\la} W_\mu$. We know $K_\la =
\cap_{0<\mu<\la} K_\mu$, hence an element $f(t)$ in
$\cap_{0<\mu<\la} W_\mu$ belongs to $K_\la$ and has its Mellin
transform vanishing at least as the zeta function does. From
$K_\la = W_\la^\prime\perp Z_\la$ it belongs to
$W_\la^\prime$. Hence $W_\la^\prime = \cap_{0<\mu<\la} W_\mu$. A
non-countable set of exceptional $\la$'s contradicts the
separability of $K$.
\end{proof}

We briefly explain how some of the considerations extend to
Dirichlet $L$-series. For an odd character the cosine transform
$\cF_+$ is replaced with the sine transform $\cF_-$, so we will
stick with an even (primitive) Dirichlet character: $\chi(-1) =
1$. Let us recall the functional equation of $L(s,\chi) =
\sum_{n\geq1} \chi(n)n^{-s}$:
$$L(s,\chi) = w_\chi q^{-s+1/2}\chi_+(s)L(1-s,\overline{\chi}) $$

where $w_\chi$ is a certain complex number of modulus $1$ and $q$
is the conductor (= period) of the primitive character
$\chi$. One has $\overline{w_\chi} = w_{\overline{\chi}}$. Tate's
Thesis \cite{Tate} gives a unified manner of deriving all these
functional equations as a corollary to the one-and-only
Poisson-Tate intertwining formula on adeles and ideles (and
additional local computations). A reference for the more
classical approach is, for example, \cite{Dav} (for easier
comparison with the classical formula, we have switched from
$L(\chi,s)$ to $L(s,\chi)$). The Poisson-Tate formula specializes
to twisted Poisson summation formulae on $\RR$, or rather on the
even functions on $\RR$ as we are dealing only with even
characters.

Let $\phi(t)$ be an even Schwartz function, and let:
$$P_\chi(\phi)(t) = \sum_{n\geq1} \chi(n)\phi(nt)$$

We suppose here that $\chi$ is not the principal character so
there is no  term $-(\int_0^\infty \phi(u)du)/|t|$ (which was
engineered to counterbalance the pole of the Riemann zeta
function at $s=1$). At the level of (right) Mellin transforms
$P_\chi$ corresponds to multiplication by $L(1-s,\chi)$.

So the composite $P_\chi\cF_+ = P_\chi\cF_+II$ acts on right
Mellin transforms as:
$$\wh\phi(s)\mapsto L(1-s,\chi)\chi_+(s)\wh\phi(1-s) = w_\chi
q^{s-1/2}L(s,\overline{\chi}) \wh\phi(1-s)$$
and this gives the $\chi$-Poisson intertwining:
$$P_\chi\cF_+  = w_\chi D_{q} IP_{\overline\chi}$$
where $D_{q}$ is the contraction of ratio $q$ which acts through
multiplication by $q^{s-1/2}$ on Mellin transforms and as
$f(t)\mapsto \sqrt{q}f(qt)$ on $L^2(0,\infty)$. Let us define the
$\chi$-co-Poisson $P_\chi^\prime$ on smooth even functions
compactly supported away from $0$ as:
$$P_\chi^\prime(\phi)(t) = \sum_{n\geq1}
\overline{\chi(n)}\frac{\phi(t/n)}{n}$$

We have $P_\chi^\prime = IP_{\overline\chi}I$, and
$P_\chi^\prime$ is the scale invariant operator with multiplier
(under the right Mellin transform) $L(s,\overline{\chi})$. From
the commutativity of $P_{\overline\chi}$ with $\Gamma_+ = \cF_+
I$ and the $\chi$-Poisson intertwining we get the
$\chi$-co-Poisson intertwining:
$$\cF_+ P_\chi^\prime = \cF_+IP_{\overline\chi}I =
P_{\overline\chi}\cF_+ = w_{\overline\chi}D_{q} IP_{\chi} =
w_{\overline\chi}D_{q} P_{\overline{\chi}}^\prime I$$

The placement of the operator $D_q$ on the right-side of the
Intertwining equation is very important\thinspace! If the even
function $\alpha(t)$ is supported in $(0,\infty)$ on $[\la_1,
\la_2]$ then its $\chi$-co-Poisson summation $f(t)$ will be
supported in $[\la_1, \infty)$ and the Fourier cosine transform
of $f(t)$ will be supported in $[1/(q\la_2),\infty)$. The product
of the lower ends of these two intervals is strictly less than
$1/q$ (if $\alpha$ is not identically zero). So this means that
we obtain (non-zero) functions which together with their cosine
transform are supported in $[\la,\infty[$ only for $\la <
1/\sqrt{q}$.

We let $W_\la^\chi \subset K_\la$ be the closure of such
$\chi$-twisted co-Poisson summations.  The Mellin transforms of
the functions in $W_\la^\chi$ are the functions
$L(s,\overline\chi)\wh{\alpha}(s)$ where $\alpha(t)$ is a smooth
function compactly supported in $[\la, \La/q]$ ($\La=1/\la$,
$\La>\sqrt{q}$). A vector $Z^\la_{\rho,k}$ is
(Hilbert-)perpendicular to $W_\la^\chi$ if and only if
$Z^\la_{\overline\rho,k}$ is Euclid-perpendicular to
$W_\la^\chi$ if and only if $\overline\rho$ is a (non-trivial)
zero of $L(s,\overline{\chi})$ of multiplicity strictly greater
than $k$, if and only $\rho$ is a (non-trivial)  zero of
$L(s,\chi)$ of multiplicity strictly greater than $k$. So:

\begin{theorem}
Let $\la< 1/\sqrt{q}$. A vector $Z^\la_{w,k}\in K_\la$ is
perpendicular to $W_\la^\chi$ if and only if $w$ is a non-trivial
zero $\rho$ of the Dirichlet $L$-function $L(s,\chi)$ of
multiplicity $m_\rho>k$.
\end{theorem}

We conclude with a statement whose analog we have already stated
and proven for the Riemann zeta function. The proof is only
slightly more involved, but as the statement is so important we
retrace the steps here.

\begin{theorem}\label{eureka}
The vectors $Z^\la_{\rho,k}\in K_\la$, $L(\rho,\chi)=0$, $0\leq
k<m_\rho$, associated with the non-trivial zeros of the Dirichlet
$L$-function (and with their multiplicities) span $K_\la$ if and
only if $\la\geq1/\sqrt{q}$.
\end{theorem}

\begin{proof}
They can not span if $\la<1/\sqrt{q}$ from the existence of
$W_\la^\chi$. Let us suppose $\la\geq1/\sqrt{q}$. Let $\wh f(s)$
be the Mellin Transform of an element of $K_\la$ which is
(Hilbert)-perpendicular to all $Z^\la_{\rho,k}$, $L(\rho,\chi) =
0$ (non-trivial), $0\leq k<m_\rho$. This says that $f(s)$
vanishes at the $\overline\rho$'s. We know already that $f(s)$
vanishes at the trivial zeros. So one has:
$$\wh f(s) = L(s,\overline\chi)\theta_1(s)$$
with an entire function $\theta_1(s)$. The image of $f$ under the
unitary $\cF_+$ will be Hilbert-perpendicular to
$\cF_+(Z^\la_{\rho,k}) = (-1)^k\,Z^\la_{1-\rho,k}$ and so
$\cF_+(f)$ is Euclid-perpendicular to the vectors associated to
the $\overline{1-\rho}$, which are the zeros of $L(s,\chi)$,
hence:
$$\wh{\cF_+(f)}(s) = L(s,\chi)\theta_2(s)$$
with an entire function $\theta_2(s)$. From  the functional
equation:
$$\wh{\cF_+(f)}(s) = \chi_+(s)\wh f(1-s)$$
we get
$$\chi_+(s)L(1-s,\overline\chi)\theta_1(1-s) =
L(s,\chi)\theta_2(s)$$
and combining with
$$L(s,\chi) = w_\chi q^{-s+1/2}\chi_+(s)L(1-s,\overline{\chi}) $$
this gives:
$$\theta_1(1-s) = w_\chi q^{-s+1/2}\theta_2(s)$$
Using Krein's theorem \cite{Kre} we deduce that $F(s) =
q^{-(s-1/2)/2}\theta_1(s)$ has finite exponential type which is
equal to
$$\max(\limsup_{\sigma\to+\infty}\frac{\log|F(\sigma)|}{\sigma},
\limsup_{\sigma\to+\infty}\frac{\log|F(1-\sigma)|}{\sigma})$$
and from $L(\sigma,\chi)\to_{\sigma\to+\infty} 1$ we see that the
exponential type of $F(s)$ is at most
$$\max(\log(\La)-\log(\sqrt{q}), \log(\La)-\log(\sqrt{q})) =
\log(\La)-\log(\sqrt{q})$$
This concludes the proof when $\la>1/\sqrt{q}$. When $\la =
1/\sqrt{q}$, we see that $q^{-(s-1/2)/2}\theta_1(s)$ has minimal
exponential type. But from $\wh f(s) =
L(s,\overline\chi)\theta_1(s)$ we deduce that $\theta_1(s)$ is
square-integrable on the line $\Real(s) = 2$. By the Paley-Wiener
theorem it thus vanishes identically.
\end{proof}

\section{Speculations on the zeta function, the renormalization group, duality}

We turn now to some speculative ideas concerning the zeta
function, the GUE hypothesis and the Riemann hypothesis. When we
wrote our (unpublished) manuscript ``The Explicit formula and a
propagator'' we had already spent some time trying to think about
the nature of the zeta function. Our conclusion, which had found
some kind of support with the conductor operator
$\log|x|+\log|y|$, stands today. The spaces $H\!P_\la$ and
especially Theorem \ref{eureka} have given us for the first time
a quite specific signal that it may hold some value. What is more
Theorem \ref{eureka} has encouraged us into trying to encompass
in our speculations the GUE hypothesis\footnote{{i.e.}  the
``Montgomery-Dyson proposal'' \cite{Mont} or ``Montgomery-Odlyzko
law'' \cite{Odly}.}, and more daring and distant yet, the Riemann
Hypothesis Herself.

We are mainly inspired by the large body of ideas associated with
the Renormalization Group, the Wilson idea of the statistical
continuum limit, and the unification it has allowed of the
physics of second-order phase transitions with the concepts of
quantum field theory. Our general philosophical outlook had been
originally deeply framed through the Niels Bohr idea of
complementarity, but this is a topic more distant yet from our
immediate goals, so we will leave this aside here.

We believe that the zeta function is analogous to a
multiplicative wave-field renormalization. We expect that there
exists some kind of a system, in some manner rather alike the
Ising models of statistical physics, but much richer in its phase
diagram, as each of the $L$-function will be associated to a
certain universality domain. That is we do not at all attempt at
realizing the zeta function as a partition function. No the zeta
function rather corresponds to some kind of symmetry
pattern\footnote{Of course in statistical physics, symmetry is
restored at high temperature and broken at low temperature. But
this is from a point of view where a continuum is considered more
symmetric than  a lattice as it has a larger symmetry group. So
here we are using the word ``symmetry'' under a more colloquial
acceptation.} appearing at low temperature. But the other
$L$-functions too may themselves be the symmetry where the system
gets frozen at low temperature.

Renormalization group trajectories flow through the entire space
encompassing all universality domains, and perhaps because there
are literally fixed points, or another more subtle mechanism,
this gives rise to sets of critical exponents associated with
each domain: the (non-trivial) zeros of the $L$-functions. So
there could be some underlying quantum dynamics, but the zeros
arise at a more classical level\footnote{``classical'' in its
kinematics: understanding the flow of coupling constants of a
quantum theory with infinitely many degrees of freedom has become
almost synonymous with understanding its ``quantum physics''.},
at the level of the renormalization group flow.

The Fourier transform as has been used constantly in this
manuscript will correspond to a simple symmetry, like exchanging
all spins up with all spins down. The functional equations
reflect this simple-minded symmetry and do not have a decisive
significance in the phase picture.

But we do believe that some sort of a much more hidden thing
exist, a Kramers-Wannier like duality exchanging the low
temperature phase with a single hot temperature phase, not
number-theoretical. If this were really the case, some universal
properties would hold across all phases, reflecting the
universality examplified by the GUE hypothesis. Of course the hot
phase is then expected to be somehow related with quantities
arising in the study of  random matrices. In the picture from
Theorem \ref{eureka}, $\la$ seems to play the r\^ole of a
temperature (inverse of coupling constant).

We expect that if such a duality did reign on our space it would
interact in such a manner with the renormalization group flow
that this would give birth to scattering processes. Indeed the
duality could be used to compare incoming to outgoing (classical)
states. Perhaps the constraints related with this interaction
would result in a property of causality equivalent to the Riemann
Hypothesis.

Concerning the duality at this time we can only picture it to be
somehow connected with the Artin reciprocity law, the ideas of
class field theory and generalizations thereof. So  here our
attempt at being a revolutionary ends in utmost conservatism.

\small

B.~Riemann, \emph{\"Uber die Anzahl der Primzahlen unter einer
gegebenen Gr\"osse}, Monatsber. Akadem. Berlin, 671-680, (1859).

\hangindent 1cm N.~Bohr, \emph{The Philosophical Writings of
Niels Bohr}, Ox Bow Press, Woodbridge, Connecticut.  \\ Volume I:
\emph{Atomic Theory and the Description of Nature} (1934)\\
Volume II: \emph{Essays 1932-1957 on Atomic Physics and Human
Knowledge} (1958)\\ Volume III: \emph{Essays 1958-1962 on Atomic
Physics and Human Knowledge} (1963)\\ Volume IV:  \emph{Causality
and Complementarity}, Supplementary papers edited by Jan Faye and
Henry J. Folse, (1999).

K.G.~Wilson, \emph{The renormalization group and critical
phenomena}, Rev. Mod. Phys. {\bf 55} (1983), 583-600.



\begin{thebibliography}{99}
\addcontentsline{toc}{section}{\numberline {8}References}

\bibitem{Bae} L.~B\'aez-Duarte,  \emph{A class of invariant
unitary operators}, Adv. in Math. {\bf 144} (1999), 1-12.

\bibitem{luisnatural} L. B\'{a}ez-Duarte, \emph{A strengthening
of the Nyman-Beurling criterion for the Riemann Hypothesis},
Rendiconti Accademia dei Lincei, to appear.

\bibitem{Bal3} {L.~B\'aez-Duarte, M.~Balazard, B.~Landreau and
E.~Saias}, \emph{Notes sur la fonction $ \zeta $ de Riemann 3},
Adv. in Math.  {\bf 149} (2000), 130-144.

\bibitem{Bal1} {M.~Balazard, E.~Saias}, \emph{The Nyman--Beurling
equivalent form for the Riemann hypothesis}, Expo. Math. {\bf 18}
(2000), 131--138.

\bibitem{Beu49} A.~Beurling, \emph{On two problems concerning
linear transformations in Hilbert space}, Acta Mathematica\ {\bf
81} (1949), 239--255.


\bibitem{Beu55} A.~Beurling, \emph{A closure problem related to
the Riemann Zeta--function}, Proc.\ Nat.\ Acad.\ Sci.\ {\bf 41}
(1955), 312-314.

\bibitem{Bra64} L.~de~Branges, \emph{Self-reciprocal functions},
J. Math. Anal. Appl. {\bf 9} (1964) 433--457.

\bibitem{Bra} L.~de~Branges, \emph{Hilbert spaces of entire
functions}, Prentice Hall Inc., Englewood Cliffs, 1968.



\bibitem{Bra86} L.~de~Branges, \emph{The Riemann hypothesis for
Hilbert spaces of entire functions},
Bull. Amer. Math. Soc. (N.S.)  {\bf 15} (1986), no. 1, 1--17.

\bibitem{Bra92} L.~de~Branges, \emph{The convergence of Euler
products}, J. Funct. Anal. 107 (1992), no. 1, 122--210.


\bibitem{Bra94} L.~de~Branges, \emph{A conjecture which implies
the Riemann hypothesis}, J. Funct. Anal. {\bf 121} (1994), no. 1,
117--184.

\bibitem{conductor} J.-F.~Burnol, \emph{The Explicit Formula and
a propagator}, 21 p., {\tt math/9809119} (1998).\\  \emph{The
Explicit Formula in simple terms}, 25 p., {\tt math/9810169}
(1998).\\ \emph{Spectral analysis of the local conductor
operator},  12 p., {\tt math/9811040} (1998).\\  \emph{Spectral
analysis of the local commutator operators}, 16 p., {\tt
math/9812012} (1998).\\  \emph{The explicit formula and the
conductor operator}, 28p. {\tt math/9902080} (1999).

\bibitem{cras1} J.-F.~Burnol, \emph{Sur les formules explicites
I: analyse invariante}, C. R. Acad. Sci. Paris {\bf 331} (2000),
S\'erie I, 423-428.

\bibitem{imrn} J.-F.~Burnol, \emph{Scattering on the p-adic field
and a trace formula}, Int. Math. Res. Not. {\bf 2000:2} (2000),
57-70.


\bibitem{jnt} J.-F.~Burnol, \emph{An adelic causality problem
related to abelian $L-$functions}, J. Number Theory {\bf 87}
(2001), no. 2, 253-269.

\bibitem{mrl} J.-F.~Burnol, \emph{Quaternionic gamma functions
and their logarithmic derivatives as spectral functions},
Math. Res. Lett. {\bf 8} (2001), no. 1-2, 209-223.

\bibitem{ams}  J.-F.~Burnol, \emph{A note on Nyman's equivalent
formulation for the Riemann Hypothesis}, in \emph{Algebraic
Methods in Probability and Statistics}, Viana M. \& Richards
D. eds, Contemp. Math. {\bf 287}, 23-26, American Mathematical
Society, Providence, RI, 2001.

\bibitem{aim}  J.-F.~Burnol, \emph{A lower bound in an
approximation problem involving the zeros of the Riemann zeta
function}, Adv. in Math. {\bf 170} (2002), 56-70.

\bibitem{cras2} J.-F.~Burnol, \emph{Sur certains espaces de
Hilbert de fonctions enti\`eres{,} li\'es \`a la transformation
de Fourier et aux fonctions L de Dirichlet et de Riemann},
C. R. Acad. Sci. Paris {\bf 333} (2001), s\'erie~I, 201-206.

\bibitem{zetaratio} J.-F.~Burnol, \emph{On an analytic estimate
in the theory of the Riemann Zeta function and a Theorem of
B\'aez-Duarte}, {\tt math/0202166}, February 2002.

\bibitem{twosys} J.-F.~Burnol, \emph{Two complete and minimal
systems associated with the zeros of the Riemann zeta function},
{\tt math/0203120}, March 2002.


\bibitem{cras3} J.-F.~Burnol, \emph{Sur les ``Espaces de Sonine''
associ\'es par de~Branges \`a la transformation de Fourier},
C. R. Acad. Sci. Paris {\bf 335} (2002), s\'erie~I, 689-692.

\bibitem{copoisson} J.-F. Burnol, \emph{Co-Poisson intertwining:
distribution and function theoretic aspects}, in preparation.

\bibitem{Con96} A.~Connes, \emph{Formule de trace en
g\'eom\'etrie non-commutative et hypoth\`ese de Riemann},
C. R. Acad. Sci. Paris {\bf 323} (1996), S\'erie I, 1231-1236.

\bibitem{Con99} A.~Connes, \emph{Trace formula in non-commutative
Geometry and the zeros of the Riemann zeta function}, Selecta
Math. (N.S.) {\bf 5} (1999) , no. 1, 29--106.

\bibitem{conreyli} J. B. Conrey, X-J. Li, \emph{A note on some
positivity conditions related to zeta and $L$-functions},
Internat. Math. Res. Notices (2000), no. 18, 929--940.

\bibitem{Dav} H.~Davenport, \emph{Multiplicative number theory},
Third edition. Revised and with a preface by Hugh L.
Montgomery. Graduate Texts in Mathematics, 74. Springer-Verlag,
New York, 2000.

\bibitem{Dym1} H.~Dym, H.P.~McKean, \emph{Fourier Series and
Integrals}, Academic Press, 1972.


\bibitem{Dym}  H.~Dym, H.P.~McKean, \emph{Gaussian processes,
function theory, and the inverse spectral problem}, Probability
and Mathematical Statistics, Vol. 31. Academic Press [Harcourt
Brace Jovanovich, Publishers], New York-London, 1976.

\bibitem{Gel} I.M.~Gel'fand, M.I.~Graev, I.I.~Piateskii-Shapiro,
\emph{Representation Theory and automorphic functions},
Philadelphia, Saunders (1969).

\bibitem{Gui}  A.P.~Guinand, \emph{A summation formula in the
theory of prime numbers}, Proc. London Math. Soc. (2) {\bf 50},
(1948), 107--119.


\bibitem{Har} S.~Haran, \emph{Riesz potentials and explicit sums
in arithmetic}, Invent. Math.  {\bf 101} (1990), 697-703.

\bibitem{Kre} M.G.~Krein, \emph{Theory of entire functions of
exponential type} (in Russian), Izv. Akad. Nauk. SSSR,
Ser. Mat. 11 (1947), No. 4, 309-326.

\bibitem{Lax59} P.~Lax, \emph{Translation invariant subspaces},
Acta Mathematica\ {\bf 101} (1959), 163--178.

\bibitem{Lax67} P.~Lax, R.~S.~Phillips, \emph{Scattering Theory},
(1st ed. 1967), Rev. Ed., Pure and Applied Mathematics, v.26,
Academic Press, 1989.

\bibitem{Mont} H.~Montgomery, \emph{The pair correlation function
of zeros of the zeta function}, Analytic Number Theory
(H.~G.~Diamond, ed.) Proc. Sympos. Pure Math., {\bf 24},
Amer. Math. Soc. (1973) 181--193.


\bibitem{Nym} B.~Nyman, \emph{On the One-Dimensional Translation
Group and Semi-Group in Certain Function Spaces}, Thesis,
University of Uppsala, 1950. 55 pp.

\bibitem{Odly} A.~M.~Odlyzko, \emph{The $10^{20}$th zero of the
Riemann zeta function and $70$ million of its neighbors}, ATT
Bell Laboratories preprint, 1989.


\bibitem{Rov66} V.~Rovnyak, \emph{Self-reciprocal functions},
Duke Math. J. {\bf 33} (1966) 363--378.

\bibitem{Rov67} J.~Rovnyak, V.~Rovnyak, \emph{Self-reciprocal
functions for the Hankel transformation of integer order}, Duke
Math. J. {\bf 34} (1967) 771--785.

\bibitem{Rov69} J.~Rovnyak, V.~Rovnyak, \emph{Sonine spaces of
entire functions}, J. Math. Anal. Appl. {\bf 27} (1969) 68--100.


\bibitem{Sonine} N.~Sonine, \emph{Recherches sur les fonctions
cylindriques et le d\'eveloppement des fonctions continues en
s\'eries}, Math. Ann. \textbf{16} (1880), 1--80.

\bibitem{Tate} J.~Tate,  \emph{Fourier Analysis in Number Fields
and Hecke's Zeta Function}, Thesis, Princeton 1950, in
\emph{Algebraic Number Theory}, {Cassels J.W.S.}, {Fr\"ohlich A.} 
eds., Academic Press, 1967.


\bibitem{TitchmarshFourier} E. C. Titchmarsh, \emph{Introduction
to the Theory of Fourier Integrals}, Clarendon Press, Oxford, 2nd
ed., 1948.

\bibitem{Titchmarsh} E. C. Titchmarsh, \emph{The Theory of the
Riemann-Zeta Function}, 2nd ed. Edited and with a preface by
D. R. Heath-Brown, Clarendon Press, Oxford, 1986.

\bibitem{Wei52} A.~Weil, \emph{Sur les ``formules explicites'' de
la th\'eorie des nombres premiers}, Comm. Sem. Math. Univ. Lund
(1952) Volume dedicated to Marcel Riesz, 252--265. Oeuvres,
Vol. II.

\bibitem{Wei72} A.~Weil, \emph{Sur les formules explicites de la
th\'eorie des nombres}, Izv. Mat. Nauk. (Ser. Mat.) {\bf 36}
(1972), 3-18. Oeuvres, Vol. III.

\bibitem{Wei74} A.~Weil, \emph{Basic Number Theory}, 3rd ed.,
Springer--Verlag, 1974.


\end{thebibliography}
\end{document}